\documentclass{article}

\usepackage{graphics}
\usepackage{amsmath,amssymb,amsthm}
\usepackage[mathscr]{eucal}

\newtheorem{theorem}{Theorem}[section]
\newtheorem{proposition}[theorem]{Proposition}

\newtheorem{lemma}[theorem]{Lemma}

\theoremstyle{definition}
\newtheorem{defi}{Definition}[section]

\newcommand{\arr}{\longrightarrow}
\newcommand{\C}{\mathbb{C}}
\newcommand{\dom}{\mathop{\mathrm{Dom}}}
\newcommand{\emp}{\varnothing}
\newcommand{\G}{\Gamma}
\newcommand{\img}[1]{\mathop{\mathrm{IMG}}\left(#1\right)}
\newcommand{\lims}[1]{\mathscr{J}_{#1}}
\newcommand{\LL}{\Lambda}
\newcommand{\M}{{\mathcal{M}}}
\newcommand{\N}{\mathbb{N}}
\newcommand{\nuke}{\mathcal{N}}
\newcommand{\prr}{\dashrightarrow}
\newcommand{\R}{\mathbb{R}}
\newcommand{\si}{\mathsf{s}}

\newcommand{\xmo}{X^{-\omega}}
\newcommand{\xs}{X^*}
\newcommand{\Z}{\mathbb{Z}}

\title{Iterated Monodromy Groups}
\author{Volodymyr Nekrashevych
\thanks{The research was supported by the Swiss
National Science Foundation and Alexander~von~Humboldt
Foundation}}

\begin{document}
\maketitle
\begin{abstract}
We associate a group $\img{f}$ to every covering $f$ of a
topological space $\M$ by its open subset.  It is the quotient of
the fundamental group $\pi_1(\M)$ by the intersection of the
kernels of its monodromy action for the iterates $f^n$. Every
iterated monodromy group comes together with a naturally defined
action on a rooted tree. We present an effective method to compute
this action and show how the dynamics of $f$ is related to the
group. In particular, the Julia set of $f$ can be reconstructed
from $\img{f}$ (from its action on the tree), if $f$ is expanding.

\end{abstract}
\section{Introduction}

The aim of this paper is to show a new connection between
dynamical systems and algebra. A group, called \emph{iterated
monodromy group} is associated to every covering $f:\M_1\arr\M$ of
a topological space by its open subset. This group encodes the
combinatorial information about the iterations of the map $f$. If
the map $f$ is expanding (\emph{hyperbolic}) then the iterated
monodromy group (together with the associated \emph{virtual
endomorphism}) contains all the ``essential'' information about
the dynamics of $f$: one can reconstruct from it the action of $f$
on its Julia set.

Let $f:\M_1\arr\M$ be a $d$-fold covering map of an arcwise
connected and locally arcwise connected topological space $\M$ by
its arcwise connected open subset $\M_1$. By $f^n$ we denote the
$n$th iteration of the mapping $f$. It is a covering of the space
$\M$ by its open subet $\M_n$. Choose an arbitrary point $t\in\M$
and let $T_t$ be the formal disjoint union of the sets $f^{-n}(t)$
of preimages of $t$ under $f^n$. The set $T_t$ has a natural
structure of a regular rooted tree with the root $t\in f^{-0}(t)$
in which every vertex $z\in f^{-n}(t)$ is connected to the vertex
$f(z)\in f^{-(n-1)}(t)$. The tree $T_t$ is called \emph{preimage
tree}.

The fundamental group $\pi_1(\M, t)$ naturally acts on each set
$f^{-n}(t)$. It is easy to see that the obtained action of
$\pi_1(\M, t)$ on $T_t$ is an action by automorphisms of the
rooted tree. This action is called \emph{iterated monodromy
action} of $\pi_1(\M)$. It does not depend, up to a conjugacy of
actions, on the choice of $t$ (Proposition~\ref{prop:img}).

The iterated monodromy action is not faithful in general.
Therefore, the following definition is introduced.

\begin{defi}
\emph{Iterated monodromy group} $\img{f}$ of the covering $f$ is
the quotient of the fundamental group $\pi_1(\M, t)$ by the kernel
of its iterated monodromy action on $T_t$.
\end{defi}

Iterated monodromy groups are discrete analogs of the following
Galois groups, defined by R.~Pink. Let $f\in\Bbbk[x]$ be a
polynomial over a field $\Bbbk$. Denote by $f^{\circ n}(x)$ its
$n$th iteration and define the polynomials $F_n(x)=f^{\circ
n}(x)-t$ over the field $\Bbbk(t)$. Let $\Omega_n$ be the
splitting field of $F_n$ over $\Bbbk(t)$ and let
$\Omega=\cup_{n\ge 1}\Omega_n$. We obtain the Galois group
$\mathop{\overline{\mathrm{IMG}}}(f)=\mathop{\mathrm{Aut}}\left(\Omega/\Bbbk(t)\right)$.

It is not hard to prove that if $f\in\C[x]$ is a post-critically
finite polynomial then the Galois group
$\mathop{\overline{\mathrm{IMG}}}(f)$ is the closure of the
iterated monodromy group $\img{f}$ in the automorphism group of
the rooted tree. Here $\img{f}$ is computed as the iterated
monodromy group of the covering $f:\M_1\arr\M$ for $\M=\C\setminus
P$, $\M_1=\C\setminus f^{-1}(P)$, where $P$ is the set of
post-critical points of $f$ (or any other finite set, for which
$f:\M_1\arr\M$ is a covering map).

We present in our paper a method to compute the action of the
iterated monodromy group $\img{f}$ on the rooted tree.

Automorphisms of rooted trees are conveniently encoded using
automata (see the survey~\cite{grineksu_en}). Any regular tree is
isomorphic to the tree $\xs$ of finite words over an alphabet $X$.
We connect two words by an edge in the tree $\xs$ if and only if
they are of the form $v$ and $vx$, where $v\in\xs$ and $x\in X$.
The root of the tree $\xs$ is the empty word $\emp$.

If $g$ is an automorphism of the rooted tree $\xs$ then for every
$x\in X$ there exists a uniquely defined automorphism $g|_x$ such
that \begin{equation} \label{eq:restriction} (xv)^g=x^gv^{g|_x}
\end{equation}
for all $v\in\xs$. This can be interpreted in terms of automata in
the following way. The automorphism $g$ is considered to be a
state of an automaton, which when reading an input letter $x\in X$
gives on output the letter $x^g$ and then changes its state to
$g|_x$. Automata of this type are called sometimes
\emph{transducers} or \emph{sequential machines} (see~\cite{eil}).

Suppose now that we have a $d$-fold covering $f:\M_1\arr\M$ of a
topological space by its open subspace. Let $t\in\M$ be a
basepoint, let $X$ be an alphabet of cardinality $d$. We choose
some bijection $\LL:X\arr f^{-1}(X)$ and paths $\ell_x$ in $\M$
connecting $t$ to $\LL(x)$.

The choice of the paths $\ell_x$ defines an isomorphism of the
rooted trees $\LL:\xs\arr T_t$, which can be used to encode the
vertices of the tree $T_t$ by words over the alphabet $X$ (see
Definition~\ref{def:naturenc}). We get in this way a
\emph{standard action} of $\img{f}$ on $\xs$, conjugating the
iterated monodromy action on $T_t$ by $\LL$ (i.e., identifying
$T_t$ with $\xs$ by $\LL$).

The standard action is computed by the following recurrent formula
(see Proposition~\ref{pr:nat}):
\[
(xv)^{\gamma}=y\left(v^{\ell_x\gamma_x \ell_y^{-1}}\right),
\]
where $\gamma_x$ is the $f$-preimage of the loop $\gamma$, which
starts at $\LL(x)$ (and ends at $\LL(y)$).

This formula can be interpreted as a description of the automaton,
whose action on $\xs$ coincides with the standard action of the
loop $\gamma$.

It is therefore possible to apply the techniques developed for the
study of groups generated by automata to the iterated monodromy
groups. We review the main definitions and results about groups
generated by automata in Subsection~\ref{ss:treexs} and in
Subsection~\ref{ss:virt} (where an algebraic description of such
actions is introduced and studied). For more details on this
topic, see the surveys~\cite{grineksu_en,bgn}.

Theory of groups generated by automata is developing intensively
in the last two decades (see the
works~\cite{grineksu_en,grigorchuk:branch,sidki_monogr,bgn} and
their bibliography). It was discovered, that many interesting
groups can be easily defined and studied using their actions on
trees.

The first example of a group of this sort was the Grigorchuk
group, defined in~\cite{grigorchuk:80_en}. It was constructed
originally as a simply defined example of an infinite finitely
generated torsion group (thus related to the General Burnside
Problem). It was discovered later that it is a group of
intermediate growth~\cite{grigorchuk:milnor_en} (and thus
answering on the Milnor's question) and that it possesses many
other interesting properties (see~\cite{grigorchuk:branch}). Later
other interesting related examples of groups acting on rooted
trees where constructed~\cite{gupta-sidki_group,bsv:jns}.

One of the main properties of these examples is the fact that the
restriction $g|_x$ (see~(\ref{eq:restriction})) is asymptotically
shorter than $g$. Such groups are called \emph{contracting}. The
contraction property provides inductive proofs of most result
about these groups. See for instance the original proof of the
fact that the Grigorchuk group is periodic
in~\cite{grigorchuk:80_en}. We discuss the basic properties of
contracting actions in Subsection~\ref{ss:contract}.

It was also discovered that groups generated by automata have rich
geometry. For example, in~\cite{gr_zu:lamp,bgr:spec} Schreier
graphs of some of such groups where described and their spectra
where computed. In particular, it became clear that the graphs of
the action of a group on the levels of the tree may converge to
some fractal space.

This observation was formalized later by the author
in~\cite{nek:lim}. It was shown that if the group action is
contracting, then a naturally defined \emph{limit space}
$\lims{G}$ together with a continuous map
$\si:\lims{G}\arr\lims{G}$ is associated to it.

The space $\lims{G}$ is defined as a quotient of the space $\xmo$
of left-infinite sequences $\ldots x_2x_1$ over the alphabet $X$
by the \emph{asymptotic equivalence relation}. Two sequences
$\ldots x_2x_1, \ldots y_2y_1\in\xmo$ are said to be
asymptotically equivalent if there exists a sequence
$\{g_k\}_{k=1}^\infty$ taking a finite number of different values
$g_k\in G$ such that
\[
(x_k\ldots x_1)^{g_k}=y_k\ldots y_1
\]
for every $k\ge 1$ (see Definition~\ref{def:asympteq}).

The asymptotic equivalence relation is described by a finite
directed labeled graph (the \emph{Moore diagram of the nucleus}):
two sequences are equivalent if and only if they are read on a
directed path of the graph (Proposition~\ref{pr:limsp}).

The shift $\sigma:\ldots x_2x_1\mapsto \ldots x_3x_2$ preserves
the asymptotic equivalence relation, hence it induces a continuous
map $\si:\lims{G}\arr\lims{G}$ on the limit space $\lims{G}$. The
obtained dynamical system $\left(\lims{G}, \si\right)$ is called
\emph{limit dynamical system} of the contracting action (see
Definition~\ref{def:limitspace}).

The main result of our paper is Theorem~\ref{th:limjul} showing
that the limit dynamical system of the iterated monodromy group of
an expanding self-covering $f:\M_1\arr\M$ is topologically
conjugate to the dynamical system $\left(\mathcal{J}(f),
f\right)$, where $\mathcal{J}(f)$ is the Julia set of $f$.

We illustrate in the last section the results of our paper on some
examples. The first class of examples are expanding endomorphisms
of Riemannian manifolds. They where studied before by M.~Schub,
J.~Franks and M.~Gromov. We show that a result of  M.~Schub and
J.~Franks, saying that an expanding endomorphism is uniquely
determined by its action on the fundamental group, is a partial
case of Theorem~\ref{th:limjul} (see Theorem~\ref{th:expman} and
Theorem~\ref{th:shub2}). This is illustrated on some concrete
examples like self-coverings of torus and Heisenberg group.
Theorem~\ref{th:limjul} and the definition of the limit space
provide an encoding of the manifolds by infinite sequences. These
encodings are interesting numeration systems on nilpotent Lie
groups (in particular on $\R^n$), which where studied by many
authors (see~\cite{cuntz_rep,vince:digtile}).

Another very interesting class of examples of iterated monodromy
groups comes from holomorphic dynamics. Every hyperbolic rational
map (i.e., a map for which orbits of all critical points are
converging to an attracting cycle) is expanding on a neighborhood
of its Julia set and Theorem~\ref{th:limjul} can be applied.

Iterated monodromy groups appeared implicitly in the
paper~\cite{lyubichminsk} of M.~Lyubich and  Y.~Minsky, since they
can be defined as the holonomy groups of the laminations, studied
in~\cite{lyubichminsk}.

The encoding of the Julia set by infinite sequences, given by
Theorem~\ref{th:limjul} was also studied before. See, for example
the papers of M.~V.~Yacobson~\cite{yacobson3,yacobson4}.

We compute several examples of iterated monodromy groups of
post-critically finite rational functions. The ``smooth'' examples
give ``usual'' groups: $\mathbb{Z}$ for $z^2$, infinite dihedral
for $z^2-2$ and a $\mathbb{Z}/2\mathbb{Z}$-extension of
$\mathbb{Z}^2$ (the group of the affine transformations $\pm x+a$,
of $\mathbb{Z}^2$) for the Latt\`es examples.

On the other hand, the group $\img{z^2-1}$ does not have a finite
presentation by defining relations (though a simple recursive
definition is known, see Theorem~\ref{th:presimgmin}). It has
exponential growth but no free non-abelian subgroups. Every its
proper quotient is solvable.

The group $\img{z^2-1}$ is the first example of an amenable group,
which can not be constructed from groups of sub-exponential growth
using the extensions and direct limits. Amenability of
$\img{z^2-1}$ was proved by L.~Bartholdi and
B.~Vir\'ag~\cite{barthvirag} using self-similarity of the random
walks on it.

It is interesting to mention that the group $\img{z^2-1}$ was
defined and studied for the first time by R.~Grigorchuk and
A.~\.Zuk (see~\cite{zukgrigorchuk:3st}) before iterated monodromy
groups where defined. They introduced this group as just an
interesting example of a group generated by a three-state
automaton. The fact that $\img{z^2-1}$ can not be obtained from
the groups of sub-exponential growth is a result
of~\cite{zukgrigorchuk:3st}.

Most results of this paper where announced in~\cite{bgn}, where
also some more examples of iterated monodromy groups where
presented.

\medskip
\noindent\textbf{Acknowledgements.} I am very grateful to
L.~Bartholdi, R.~Grigorchuk and P.~de~la~Harpe for discussions,
help and input into this work.

I would also like to thank R.~Pink for productive discussions
which eventually led to the definition of iterated monodromy
groups.

\section{Preliminary definitions}
\label{s:prelim}

\subsection{Rooted trees} \label{ss:rooted}

A \emph{rooted tree} $T$ is a simplicial graph without cycles with
a marked vertex $v_\emp$ called the \emph{root}. An
\emph{isomorphism} of two rooted trees is an isomorphism of the
graphs, which preserves the roots.

Every rooted tree $T$ can be defined by a sequence of sets and
maps
\begin{equation}
\label{eq:invseq}
X_0\stackrel{f_1}{\longleftarrow}X_1\stackrel{f_2}{\longleftarrow}X_2\cdots,
\end{equation}
where $X_0=\{v_\emp\}$ contains only the root. Here $V=\cup_{n\ge
0} X_n$ is the set of vertices of $T$ and every vertex $v\in X_n$
is connected by an edge with the vertex $f_n(v)\in X_{n-1}$.

The set $X_n$ is called the $n$th \emph{level} of the tree $T$ and
is uniquely defined as the set of the vertices which are on
distance $n$ from the root.

The rooted $T$ tree is \emph{$d$-regular} if every point $v\in
X_n$ has exactly $d$ preimages under $f_{n+1}$ for every $n$.

Let $X$ be a finite set (an \emph{alphabet}). Denote by $\xs$ the
free monoid generated by $X$, i.e., the set of all finite words of
the form $x_1x_2\ldots x_n$, where $x_i\in X$, (including the
empty word $\emp$). It has a natural structure of a tree, where
every word $v\in\xs$ is connected to the word $vx$ for every $x\in
X$. In this way we get a rooted tree with the root in $\emp$. It
is the tree defined by the sequence
$$
X^0\stackrel{f_1}{\longleftarrow}X^1\stackrel{f_2}{\longleftarrow}X^2\cdots,
$$
where $f_n(x_1x_2\ldots x_n)=x_1x_2\ldots x_{n-1}$ (and
$X^0=\{\emp\}$ consists only of the empty word).

The $n$th level of the tree $\xs$ is the set $X^n$ of the words of
\emph{length $n$}. We denote the length of a word $v\in\xs$ by
$|v|$, so that $v\in X^n$ if and only if $|v|=n$.

The tree $\xs$ is $d$-regular for $d=|X|$ and every $d$-regular
tree is isomorphic to $\xs$.

It is easy to see that automorphisms of rooted trees preserve the
levels and that the following simple lemma holds.

\begin{lemma}
\label{l:isom} A bijection $g:\cup_{n\ge 0}X_n\arr\cup_{n\ge
0}X_n$ is an automorphism of the rooted tree defined by the
inverse sequence~(\ref{eq:invseq}), if and only if $(X_n)^g=X_n$
for every $n\ge 0$ and
$$
(f_n(v))^g=f_n\left(v^g\right)
$$
for every $v\in X_n$.\qed
\end{lemma}

\subsection{Partial self-coverings}

We will use the standard terminology and facts about the covering
maps (see, for example~\cite{massey:algtop}).

\begin{defi}
\label{def:psc} Let $\M$ be an arcwise connected and locally
arcwise connected topological space. A $d$-fold \emph{partial
self-covering map} on the space $\M$ is a $d$-fold covering map
$f:\M_1\arr\M$, where $\M_1$ is an open arcwise connected subset
of $\M$.
\end{defi}

Recall that $f:\mathcal{M}_1\arr\mathcal{M}_2$ is a $d$-fold
covering map if it is surjective and every point
$x\in\mathcal{M}_2$ has a neighborhood $\mathcal{U}$ such that the
preimage $f^{-1}(\mathcal{U})$ is a disjoint union of $d$ subsets
$\mathcal{U}_i$ for which the restriction
$f:\mathcal{U}_i\arr\mathcal{U}$ is a homeomorphism.

\medskip\noindent\textbf{Examples.}
{\bf 1. Self-covering} is a covering $f:\M\arr\M$ of a space by
itself. As a simple example, consider the double self-covering
$x\mapsto 2x$ of the circle $\R/\Z$ (equivalently, the map
$z\mapsto z^2$ on the unit circle $\{z\in\C\;:\;|z|=1\}$).

\noindent {\bf 2. Branched coverings.} Let $\widehat{\M}$ be a
topological space. A map $f:\widehat{\M}\arr\widehat{\M}$ is a
\emph{branched covering} if there exists a set
$R\subset\widehat{\M}$ of \emph{branching points} such that $f$ is
a local homeomorphism in every point $x\in\widehat{\M}\setminus
R$. Then the set $P=\cup_{k=0}^\infty f^k(R)$ is called the
\emph{postcritial set}. If the set
$\M=\widehat{\M}\setminus\overline{P}$ is arcwise connected and
locally arcwise connected, then $f:\M_1\arr\M$ is a partial
self-covering of the set $\M$, where $\M_1=f^{-1}(\M)$.  Here
$\overline{P}$ denotes the closure of the set $P$.

For example, the famous theorem of Thurston considers
\emph{postcritically finite} branched coverings
$f:\hat\C\arr\hat\C$ of the complex sphere, i.e., the branched
coverings for which the set $P$ is finite
(see~\cite{DH:Thurston}).

\noindent {\bf 3. Rational functions.} In particular, a rational
function $f\in\C(z)$ defines a branched covering of the complex
sphere $\hat\C=\C\cup\infty$. The set of branching points of a
rational function is the set of its \emph{critical values}, i.e.,
the values of the function $f$ in the critical points of the
function. If the postcritical set $P$ is small enough (for
instance, if it is finite or has a finite number of accumulation
points), then the function $f$ is a partial self-covering of the
set $\M=\hat\C\setminus\overline{P}$.

\noindent {\bf 4. Polynomial-like maps.} Let $U$ and $V$ be open
disks in $\hat{\C}=\C\cup\{\infty\}$. A holomorphic map $f:U\arr
V$ is said to be \emph{proper}, if an $f$-preimage of every
compact subset $K\subset V$ is a compact subset of $U$.

The following notion was introduced in~\cite{dh:plike}.
\begin{defi}
 A \emph{polynomial-like map} $f:U\arr V$ is a proper map between open
disks such that the closure of $U$ is a compact subset of $V$.
\end{defi}

Let $\overline P$ be the closure of the set of post-critical
points of the polynomial-like map and suppose that the sets
$V\setminus\overline P$ and $f^{-1}\left(V\setminus\overline
P\right)$ are arcwise connected. Then the polynomial-like map $f$
is a partial self-covering of the set $\M=V\setminus\overline P$.

\section{The definition}

\subsection{Iterated monodromy groups}

Let us fix some $d$-fold partial self-covering $f$ of an arcwise
connected and locally arcwise connected space $\M$. We have the
following classical
\begin{lemma}
\label{l:paths} For every path $\gamma$ in $\M$ and every
$f$-preimage $z$ of the beginning of $\gamma$ there exists a
unique path $\gamma'$ in $\M_1$ beginning at $z$ and such that
$f(\gamma')=\gamma$.
\end{lemma}

\noindent \textbf{Notation.} If $\gamma$ is a path in $\M$ and a
point $z\in\M$ is such that $f^n(z)$ is the beginning of $\gamma$,
then we denote by $f^{-n}(\gamma)[z]$ the preimage of $\gamma$
under $f^n$, which starts at the point $z$.
\medskip

Let $t\in \M$ be an arbitrary point. We get an inverse sequence
$$
\{t\}\stackrel{f}{\longleftarrow}
f^{-1}(t)\stackrel{f}{\longleftarrow} f^{-2}(t)
\stackrel{f}{\longleftarrow}
f^{-3}(t)\stackrel{f}{\longleftarrow}\cdots
$$
defining a rooted tree $T_t$ called the \emph{preimages tree} of
the point $t$. The tree $T_t$ is $d$-regular, since the map $f$ is
a $d$-fold covering.

Let now $\gamma$ be a loop in $\M$ based at $t$, i.e., a path
starting and ending at $t$.  For every vertex $z\in f^{-n}(t)$ of
the $n$th level of the preimage tree $T_t$  denote by $z^\gamma$
the end of the path $f^{-n}(\gamma)[z]$.

Then we obviously have $f^n\left(z^\gamma\right)=t$, so the
element $z^\gamma$ also belongs to the $n$th level of $T_t$.

\begin{proposition}
  \label{prop:img}
The map $z\mapsto z^\gamma$ is an automorphism of the preimage
tree, which depends only on the homotopy class of $\gamma$ in
$\M$. In this way we get an action of the fundamental group
$\pi_1(\M, t)$ on the tree $T_t$. Up to a conjugacy,  the action
does not depend on the choice of the basepoint.
\end{proposition}

We say that an action of a group $G_1$ on a set $M_1$ is conjugate
to an action of a group $G_2$ on a set $M_2$ if there exists an
isomorphism $\phi:G_1\arr G_2$ and a bijection $l:M_1\arr M_2$
such that
\begin{equation}
\label{eq:isomact} l\left(x^g\right)=l\left(x\right)^{\phi(g)}
\end{equation}
for every $x\in M_1$ and $g\in G$.

\begin{proof}
The fact that the map $z\mapsto z^\gamma$ defines for every $n\ge
1$ an action of the fundamental group $\pi(\M, t)$ on the set
$f^{-n}(t)$ is classical (see, for example, \S\ 7
of~\cite{massey:algtop}).

Let us prove that the map $g:z\mapsto z^\gamma$ is an automorphism
of the tree $T_t$, using Lemma~\ref{l:isom}. The map $g$ defines a
permutation of every level $f^{-n}(t)$ of the tree $T_t$. Let
$z\in f^{-n}(t)$ be an arbitrary point of the $n$th level. The
point $z^\gamma$ is the end of the path
$\gamma'=f^{-n}(\gamma)[z]$. Then $f(\gamma')$ is equal to
$f^{-(n-1)}(\gamma)[f(z)]$, thus the end of $f(\gamma')$ is equal
to $f(z)^\gamma$. But it is also obviously equal to
$f\left(z^\gamma\right)$. Therefore,
$$
f\left(z\right)^\gamma=f\left(z^\gamma\right)
$$
for every $z\in f^{-n}(t)$ and Lemma~\ref{l:isom} shows that the
map $z\mapsto z^\gamma$ is an automorphism of the tree $T_t$.

Let $t'$ be another basepoint. Choose a path $\ell$ starting at
$t$ and ending at $t'$. Let us define a map $l: T_t\arr T_{t'}$,
which maps every point $z\in f^{-n}(t)$ to the end of the path
$f^{-n}(\ell)[z]$. Same considerations as above show that $l$ is
an isomorphism of the rooted trees. The path $\ell$ defines also
an isomorphism $\phi:\pi(\M, t)\arr\pi(\M, t')$ of the fundamental
groups by the formula $\phi(\gamma)=\ell^{-1}\gamma\ell$.

Let $\gamma$ be an arbitrary loop at $t$ and let $z\in f^{-n}(t)$
be an arbitrary vertex of the $n$th level of the tree $T_t$. Then
$z^\gamma$ is the end of the path $f^{-n}(\gamma)[z]$.  The path
$f^{-n}(\ell)\left[z^\gamma\right]$ begins at $z^\gamma$ and ends
in $l\left(z^\gamma\right)$. The path $f^{-n}(\ell)[z]$ begins at
$z$ and ends in $l(z)$.

Consequently, $\left(f^{-n}(\ell)[z]\right)^{-1}\cdot
f^{-n}(\gamma)[z]\cdot f^{-n}(\ell)\left[z^\gamma\right]$ is a
path, starting at $l(z)$, ending at $l\left(z^\gamma\right)$ and
equal to
$f^{-n}\left(\ell^{-1}\gamma\ell\right)\left[l(z)\right]$. This
implies that $l\left(z^\gamma\right)=l(z)^{\phi(\gamma)}$.\qed
\end{proof}

\begin{defi}
  \label{defi:img}
The action of the fundamental group $\pi_1(\M, t)$ on the preimage
tree $T_t$, described in Proposition~\ref{prop:img}, is called
\emph{iterated monodromy action} for the map $f$. The quotient of
$\pi_1(\M, t)$ by the kernel of this action is called the
\emph{iterated monodromy group (i.m.g.)} of the map $f$, denoted
$\img{f}$.
\end{defi}

Proposition~\ref{prop:img} implies that the group $\img{f}$ and
its action on the tree do not depend on the choice of the point
$t$.


\subsection{Standard actions of iterated monodromy groups on $\xs$}

If we want to compute the action of the iterated monodromy group
on $T_t$, then we have to find some convenient way to encode the
vertices of the tree $T_t$ as finite words over an alphabet $X$ of
$d$ letters. We show here a class of naturally defined enconding
which use the paths in the set $\M$. With respect to these
encodings the iterated monodromy action will be self-similar,
which will make the methods, developed for self-similar groups of
automata, applicable to the iterated monodromy groups. We will
discuss self-similar actions in general later.

Let us consider an alphabet $X$ with $d$
letters together with a bijection $\LL:X\arr f^{-1}(t)$. For every
$x\in X$ we choose a path $\ell_x$ in $\M$, starting at $t$ and
ending at $\LL(x)$.

\begin{defi}\label{def:naturenc}
The isomorphism $\LL: \xs\arr T_t$ is defined putting
$\LL(\emp)=t$, and then inductively putting $\LL(xv)$ to be equal
to the end of the path
$$f^{-(n-1)}(\ell_x)\left[\LL(v)\right],$$ where $v\in X^{n-1}$
and $x\in X$.
\end{defi}

Note that $f^{n-1}\left(\LL(xv)\right)$ is the end of the path
$\ell_x$, so that $f^n(z)=t$ and $\LL(xv)$ belongs to the $n$th
level of the preimage tree.

\begin{proposition}
\label{pr:isomorphism} The constructed map $\LL:\xs\arr T_t$ is an
isomorphism of the rooted trees.
\end{proposition}

\begin{proof}
It follows from the construction that the map $\LL$ preserves the
levels of the trees and is surjective on them. It is consequently
a bijection. Let us prove by induction on $n$ that the equality
$$
f\left(\LL(vx)\right)=\LL(v)
$$
holds for all $v\in X^n$ and $x\in X$. This will imply that the
map $\LL$ preserves the vertex adjacency and thus is an
isomorphism.

The equality is true for $n=0$. Suppose that it holds for $n=k$.
Let $v\in X^k$ and $x, y\in X$ be arbitrary. We are going to prove
that $f(\LL(yvx))=\LL(yv)$. The end of the path
$\gamma_1=f^{-(n-1)}(\ell_y)\left[\LL(v)\right]$ is, by
definition, $\LL(yv)$. The end of the path
$\gamma_2=f^{-n}(\ell_y)\left[\LL(vx)\right]$ is $\LL(yvx)$. By
assumption, $f(\LL(vx))=\LL(v)$, so we get $f(\gamma_2)=\gamma_1$,
hence, for their endpoints: $f(\LL(yvx))=\LL(yv)$. \qed\end{proof}

\begin{defi}
\label{d:nat}
 The \emph{standard action} of the group $\img{f}$ (or of $\pi_1(\M, t)$)
on the tree $\xs$ is the action obtained from the action on the
preimage tree $T_t$ conjugating it by the isomorphism $\LL:\xs\arr
T_t$, i.e., the action
\[
v^g=\LL^{-1}\left(\LL(v)^g\right).
\]
\end{defi}

The standard actions can be computed using the following recurrent
formulae.

\begin{proposition}
\label{pr:nat} Let $L=\{\ell_x\}$ be a collection of paths
defining a standard action of $\pi_1(\M, t)$ on $\xs$. Then we
have for $\gamma\in\pi_1(\M, t), x\in X, v\in\xs$ the following
relation
\begin{equation}
\label{eq:nat} (xv)^{\gamma}=y\left(v^{\ell_x\gamma_x
\ell_y^{-1}}\right),
\end{equation}
where $\gamma_x=f^{-1}(\gamma)\left[\LL(x)\right]$ and
$y=x^\gamma$.
\end{proposition}

%

\begin{proof}
Note that the path $\ell_x\gamma_x\ell_y^{-1}$ is obviously a loop
based at $t$.

Let $|v|=n$ and $\LL(yu)=\left(\LL(xv)\right)^\gamma$, where $y\in
X$ and $u\in X^n$. Denote
$$
\ell^v_{xv}=f^{-n}(\ell_x)\left[\LL(v)\right],\quad
\ell^u_{yu}=f^{-n}(\ell_y)\left[\LL(u)\right],
$$
and let
$$
\gamma_{xv}=f^{-(n+1)}(\gamma)\left[\LL(xv)\right].
$$
Then the end of $\ell^v_{xv}$ is equal to $xv$ and the end of
$\ell^u_{yu}$ is equal to $yu$. The end of $\gamma_{xv}$ is also
equal to $yu$, thus we get a path
$\ell^v_{xv}\gamma_{xv}\left(\ell^u_{yu}\right)^{-1}$ from
$\LL(v)$ to $\LL(u)$. Its image under $f^n$ is the path
$\ell_x\gamma_x\ell_y^{-1}$, thus
$u=v^{\ell_x\gamma_x\ell_y^{-1}}$. \qed\end{proof}

\subsection{The tree $\xs$ and automata}
\label{ss:treexs}

We recall here some basic facts about automorphisms of the rooted
tree $\xs$, Meely automata (transducers) and self-similar groups.
More details can be found
in~\cite{grigorchuk:branch,grineksu_en,sidki_monogr,sush:avt_en,bgn}.

Lemma~\ref{l:isom} states that a map $g:\xs\arr\xs$ is an
automorphism if and only if it preserves the length of the words
and for every $x_1x_2\ldots x_n\in X^n$ there exists $y_n\in X$
such that $(x_1x_2\ldots x_n)^g=(x_1x_2\ldots x_{n-1})^gy_n$. It
follows that for every $v\in\xs$ and $u\in X^n$ there exists $w\in
X^n$ such that $(vu)^g=v^gw$. Moreover, Lemma~\ref{l:isom} implies
that the map $u\mapsto w$ is again an automorphism of the rooted
tree $\xs$. We denote this automorphism by $g|_v$ and call it
\emph{restriction} of $g$ at $v$.

The following properties of restrictions hold:
\begin{eqnarray}
\label{eq:restr1}
(vu)^g &=& v^gu^{g|_v}\\
\label{eq:restr2}
g|_{vu} &=& \left(g|_v\right)|_u\\
\label{eq:restr3}
(g_1g_2)|_v &=& \left(g_1|_v\right)\left(g_2|_{v^{g_1}}\right)\\
\label{eq:restr4} \left(g^{-1}\right)|_v &=&
\left(g|_{v^{g^{-1}}}\right)^{-1}.
\end{eqnarray}

The set $x_1\xs$ of all the words starting with a fixed letter
$x_1\in X$ is a subtree of $\xs$ defined by the inverse sequence
$$
\{x_1\}\stackrel{f_2}{\longleftarrow}x_1X^1\stackrel{f_3}{\longleftarrow}x_1X^2\cdots.
$$

Every automorphism $g$ of the rooted tree $\xs$ acts on every
subtree $x\xs$ by the automorphism $g|_x$ and then permutes the
subtrees $x\xs$ by the permutation induced by $g$ on the set
$X^1\subset\xs$.

This leads to an interpretation of the automorphisms of rooted
trees in terms of automata.

\begin{defi}
An \emph{automaton} over the alphabet $X$ is a triple
$\left\langle Q, \lambda, \pi\right\rangle$, where
  \begin{enumerate}
  \item $Q$ is a set (\emph{the set of the internal states});
  \item $\lambda:Q\times X\arr X$ is a map, called \emph{the output function};
  \item $\pi:Q\times X\arr Q$ is a map, called \emph{the transition function}.
  \end{enumerate}

  An automaton is \emph{finite} if the set $Q$ is finite.
\end{defi}

It is convenient to define automata by their \emph{Moore
diagrams}.

\begin{defi}
A \emph{Moore diagram} of an automaton $A=\left\langle Q, \lambda,
\pi\right\rangle$ is a labeled directed graph with the set of
vertices $Q$ and the set of arrows $Q\times X$, where $(q, x)$ is
an arrow starting at $q$, ending at $\pi(q, x)$, and labeled by
the pair $(x, \lambda(q, x))\in X\times X$.
\end{defi}

On Figure~\ref{fig:nuc} an example of a Moore diagram is shown.

We interpret an automaton $A=\left\langle Q, \lambda,
\pi\right\rangle$ as a machine, which being in a state $q\in Q$
and reading on the input tape a letter $x$ goes to the state
$\pi(q, x)$ and prints on the output tape the letter $\lambda(q,
x)$. It processes the words in this way, so that the automaton $A$
with the initial state $q$ defines a map $A_q:\xs\arr\xs$ by
inductive formula
$$
A_q(\emp)=\emp,\quad A_q(xv)=\lambda(q, x) A_{\pi(q, x)}(v),
$$
where $x\in X$ and $v\in \xs$ are arbitrary.

It follows from Lemma~\ref{l:isom} that if the map $A_q$ is
bijective, then it is an automorphism of the tree $\xs$.

On the other hand~\ref{eq:restr1} implies that if $g$ is an
automorphism of the tree $\xs$, then it is defined by the
automaton $A=\left\langle Q, \lambda, \pi\right\rangle$ with the
initial state $g=g|_\emp$, where $Q=\{g|_v: v\in\xs\}$ is the set
of all possible restrictions of $g$ and the maps $\pi$ and
$\lambda$ are defined by
$$
\pi(g, x)=g|_x, \quad \lambda(g, x)=x^g.
$$

An automorphism $g$ of the tree $\xs$ is said to be \emph{finite
state} if it is defined by a finite automaton, i.e., if the set
$\{g|_v: v\in\xs\}$ is finite. The set $\mathcal{FA}$ of all
finite state automorphisms is a countable subgroup of the
automorphism group $\mathcal{A}$ of the rooted tree $\xs$. The
group $\mathcal{FA}$ is called the \emph{group of finite
automata}.

An important class of automorphism groups of the tree $\xs$ is the
class of \emph{self-similar}, or \emph{state-closed} groups.

\begin{defi}
\label{def:ssim} An action of a group $G$ on the tree $\xs$ is
said to be \emph{self-similar} if for every $g\in G$ and $x\in X$
there exists $h\in G$ and $y\in X$ such that
$$
(xw)^g=y\left(w^h\right)
$$
for every $w\in \xs$.
\end{defi}

An automorphism group $G$ of the tree $\xs$ is self-similar if and
only if for all $g\in G$ and $x\in X$ the restriction $g|_x$
belongs to $G$ (from this originates the often used term
``\emph{state-closed group}'').

Self-similar groups can be also defined as the groups, generated
by automata, in the sense of the following definition (used for
the first time in~\cite{grigorchuk:semigr_en}).

\begin{defi}
\label{def:genaut} Let $A=\left\langle Q, \lambda,
\pi\right\rangle$ be an automaton, such that the transformation
$A_q$ is invertible for every $q\in G$. Then the automorphism
group of the tree $\xs$ generated by the transformations $A_q$,
$q\in Q$ is called the \emph{group, generated by the automaton
$A$} and is denoted $\langle A\rangle$.
\end{defi}

Another interpretation of self-similar actions uses the following
notion of a permutational wreath product.

\begin{defi}
Let $H$ be a group acting by permutations on a set $X$ and let $G$
be an arbitrary group. Then the \emph{permutational wreath
product} $G\wr H$ is the semi-direct product $G^X\rtimes H$, where
$H$ acts on the direct power $G^X$ by the respective permutations
of the direct factors.
\end{defi}

The elements of a permutational wreath product $G\wr H$ are
written in the form $g\cdot h$, where $h\in H$ is a permutation of
$X$ and $g\in G^X$ is a function from $X$ to $G$. If we fix some
indexing $X=\{x_1, x_2, \ldots, x_d\}$ of the set $X$, then the
elements $g\in G^X$ can be written as tuples $g=(g_1, g_2, \ldots,
g_d)$, where $g_i=g(x_i)$.

For every self-similar automorphism group $G$ of the tree $\xs$ we
get a naturally defined homomorphism $\psi:G\arr G\wr S(X)$, where
$S(X)$ is the symmetric permutation group of the alphabet $X$.
This homomorphism is defined by the formula $\psi(g)=\tilde g\cdot
\alpha_g$, where $\alpha_g\in S(X)$ is the restriction of the
automorphism $g$ onto the set $X=X^1\subset\xs$, and $\tilde g\in
G^X$ is defined as $\tilde g(x)=g|_x$, or, in other words
$$
\psi(g)=(g|_{x_1}, g|_{x_2}, \ldots, g|_{x_d})\alpha_g.
$$

The homomorphism $\psi$ is called the \emph{permutational
recursion}, associated to the self-similar action. Note that the
permutational recursion is injective. In particular, we have an
isomorphism $\mathcal{A}\cong\mathcal{A}\wr S(X)$ for the
automorphism group $\mathcal{A}$ of the tree $\xs$. We will
usually identify the elements $g\in\mathcal{A}$ with their images
in $\mathcal{A}\wr S(X)$ under the permutational recursion $\psi$.

In the other direction, if we have a homomorphism $\psi:G\arr G\wr
S(X)$, then it defines a self-similar action of the group $G$ on
the tree $\xs$ by the recurrent formula
$$
(xv)^g=x^{\alpha_g}v^{\tilde g(x)},\quad \emp^g=\emp,
$$
where $\alpha_g\in S(X)$ is a permutation and $\tilde g\in G^X$ is
a function $X\arr G$ such that $\psi(g)=\tilde g\cdot \alpha_g$.
Note, that this self-similar action is not faithful in general.

\subsection{Virtual endomorphisms}
\label{ss:virt}

Here we recall some facts about virtual endomorphisms of groups.
For more details see the
papers~\cite{nek:stab,nek:lim,nek:iterat}.

Every partial self-covering $f$ of a topological space induces a
virtual endomorphism $\phi_f$ of the fundamental group of the
space. The dynamics of the self-covering is very closely related
to the dynamics of $\phi_f$.

\begin{defi}
A \emph{virtual homomorphism} $\phi:G_1\prr G_2$ is a homomorphism
$\phi:\dom\phi\arr G_2$, where $\dom\phi\le G_1$ is a subgroup of
finite index, called the \emph{domain} of the virtual
homomorphism. A \emph{virtual endomorphism} of a group $G$ is a
virtual homomorphism $\phi:G\prr G$.
\end{defi}

It is not hard to prove that a composition of two virtual
homomorphisms is again a virtual endomorphism.

\paragraph{Example 1.}
Let us take a group $G$ with a faithful self-similar action on the
tree $\xs$ (see Definition~\ref{def:ssim}). Suppose that the
action is transitive on the first level $X^1$. Then for any
$x_0\in X$ we have the \emph{associated virtual endomorphism} of
the action, defined as
$$
\phi_{x_0}(g)=g|_{x_0}
$$
where $\dom\phi_{x_0}$ is equal to the stabilizer $G_{x_0}$ of the
element $x_0$. The subgroup $G_{x_0}$ has index $|X|$ in $G$. The
fact that the map $\phi_{x_0}$ is a virtual endomorphism follows
from (\ref{eq:restr3}).

\paragraph{Example 2.}
Let $f:\M_1\arr\M$ be a $d$-fold partial self-covering map of the
space $\M$ and let $\ell_x$ be a path, connecting the basepoint
$t$ with one of its $f$-preimages $x$. The set
$G_1\subset\pi_1(\M, t)$ of all loops $\gamma$ such that the path
$f^{-1}(\gamma)[x]$ is again a loop, is a subgroup of index $d$ in
$\pi_1(\M, t)$, which is isomorphic to the fundamental group
$\pi_1(\M_1, x)$. Let $\phi_f:\pi_1(\M, t)\prr\pi_1(\M, t)$ be the
virtual endomorphism, equal to the composition of the
homomorphisms
$$
\pi_1(\M, t)>G_1\longrightarrow\pi_1(\M_1,
x)\longrightarrow\pi_1(\M, x)\longrightarrow\pi_1(\M, t),
$$
where
\begin{itemize}
\item[(i)]
$G_1\arr\pi_1(\M_1, x)$ is the isomorphism, which carries a loop
to its $f$-preimage,
\item[(ii)]
$\pi_1(\M_1, x)\arr\pi_1(\M, x)$ is the homomorphism induced by
the inclusion $\M_1\subseteq \M$,
\item[(iii)]
$\pi_1(\M, x)\arr\pi_1(\M, t)$ is the isomorphism defined by the
path $\ell_x$, i.e., the map $\gamma\mapsto \ell_x\gamma
\ell_x^{-1}$.
\end{itemize}

We say that the virtual endomorphism $\phi_f$ of the fundamental
group $\pi_1(\M, t)$ is \emph{induced} by the map $f$ and the path
$\ell_x$.
\medskip

\begin{defi}
Two virtual endomorphisms $\phi_1, \phi_2$ of a group $G$ are said
to be \emph{conjugate} if there exist $g, h\in G$ such that
$\dom\phi_1=g^{-1}\cdot\dom\phi_2\cdot g$ and
$\phi_2(x)=h^{-1}\phi_1(g^{-1}xg)h$ for every $x\in\dom\phi_2$.
\end{defi}

In Example 1, if we take two different letters $x,y \in X$ and
define the respective associated virtual endomorphisms $\phi_x$
and $\phi_y$, then they will be conjugate, since there will exist
an element $h\in G$ such that $x^h=y$, and then
$$
(yv)^{h^{-1}gh}=y\left(v^{h|_x^{-1}g|_xh|_x}\right)
$$
for every $g\in\dom\phi_x$ and every $v\in\xs$. Thus,
$$\phi_x(g)=h|_x\phi_y(h^{-1}gh)h|_x^{-1},$$ and the virtual
endomorphisms $\phi_x$ and $\phi_y$ are conjugate.

For the case of Example 2 we have the following proposition.

\begin{proposition}
If the virtual endomorphisms $\phi_1$ and $\phi_2$ are induced by
the map $f$ and the paths $\ell_1$ and $\ell_2$ respectively, then
they are conjugate.
\end{proposition}

\begin{proof}
Let $x_1$ be the end of the path $\ell_1$ and $x_2$ be the end of
the path $\ell_2$. Both points $x_1$ and $x_2$ belong to
$f^{-1}(t)$. There exists a path $\rho$ starting at $x_2$ and
ending at $x_1$. Its image $h=f(\rho)$ is a loop at $t$. Let
$\gamma$ be an arbitrary element of $\dom\phi_2$, i.e., such a
loop at $t$, that its preimage $\gamma_{x_2}=f^{-1}(\gamma)[x_2]$
is a loop at $x_2$.

We have:
\begin{eqnarray*}
\phi_2(\gamma) & = & \ell_2\gamma_{x_2}\ell_2^{-1}\\
\phi_1(h^{-1}\gamma h) & = &
\ell_1\rho^{-1}\gamma_{x_2}\rho\ell_1^{-1}.
\end{eqnarray*}

Therefore,
$$
\phi_2(\gamma)=\left(\ell_2\ell_1^{-1}\rho\right)\phi_2(\gamma)\left(\ell_1\rho^{-1}\ell_2^{-1}\right)=
g^{-1}\phi_1(\gamma) g,
$$
where $g=\ell_1\rho^{-1}\ell_2^{-1}$ is a loop at $t$.
\qed\end{proof}

%

Let us show how to reconstruct the self-similar action from the
associated virtual endomorphism.

Recall, that if $H$ is a subgroup of a group $G$, then a
\emph{right coset transversal} for $H$ is a set $T\subset G$ such
that $G$ is a disjoint union of the sets $Hg$, $g\in T$.

\begin{defi}
\label{def:ssf} Let $\phi$ be a virtual endomorphism of a group
$G$, let $T$ be a right coset transversal $\{r_x\}_{x\in X}$ for
the subgroup $\dom\phi$ and let $C$ be a sequence $\{h_x\}_{x\in
X}$ of elements of the group $G$. An \emph{action defined by the
triple} $(\phi, T, C)$ is the action of the group $G$ by
automorphisms of the rooted tree $\xs$ defined recurrently by the
formula
\begin{equation}
\label{eq:nat2} (xw)^g=y w^{h_x^{-1}\phi(r_xgr_y^{-1})h_y},
\end{equation}
where $x\in X$, $w\in\xs$, $g\in G$ and the element $y\in X$ is
defined by the condition $r_xgr_y^{-1}\in\dom\phi$.

An action \emph{defined by the pair} $(\phi, T)$ is the action,
defined by the triple $(\phi,T, C_0)$ where $C_0=\{1, 1, \ldots
1\}$.
\end{defi}

Note, that if $\phi$ is onto, then the action, defined by the
triple $(\phi, T, C)$ is defined by the pair $(\phi, T')$, where
$T'=\{r_x'=g_x^{-1}r_x\}_{x\in X}$, for $g_x$ such that
$\phi(g_x)=h_x$.



\begin{proposition}
Suppose that $\phi$ is a virtual endomorphism of a group $G$,
$T=\{r_x\}_{x\in X}$ is a right coset transversal for the subgroup
$\dom\phi$ and $C=\{h_x\}_{x\in X}$ is a sequence of elements of
$G$. Then (\ref{eq:nat2}) gives is a well defined action of the
group $G$ on $\xs$.
\end{proposition}

\begin{proof}
Formula~(\ref{eq:nat2}) gives a well defined action of $G$ on the
first level $X^1$ of the tree $\xs$, conjugate to the natural
action of $G$ on the right cosets $\dom\phi\cdot h$. Let $g_1,
g_2\in G$ be two arbitrary elements. Suppose that $r_xg_1r_y^{-1}$
and $r_yg_2r_z^{-1}$ belong to $\dom\phi$. Then
$r_xg_1g_2r_z^{-1}$ also belongs to $\dom\phi$ and
\[
h_x^{-1}\phi(r_xg_1g_2r_z^{-1})h_z=h_x^{-1}\phi(r_xg_1r_y^{-1})h_y\cdot
h_y^{-1}\phi(r_yg_2r_z^{-1})h_z.
\]
It follows now by induction on $n$ that formula~(\ref{eq:nat2})
gives a well defined action of the group $G$ on the $n$th level
$X^n$ of the tree $\xs$. \qed\end{proof}

\begin{proposition}
Every faithful self-similar action is defined by the associated
virtual endomorphism $\phi_{x_0}$ , a coset transversal
$T=\{r_x\}_{x\in X}$ such that $(x_0)^{r_x}=x$ and the sequence
$C=\{h_x=r_x|_{x_0}\}_{x\in X}$.
\end{proposition}

\begin{proof}
Take any $x\in X$, $g\in G$ and suppose that $y\in X$ and $h\in G$
are such that $(xw)^g=yw^h$ for all $w\in\xs$. Such $y$ and $h$
are uniquely defined, namely $y=g(x)$ and $h=g|_x$. Then:
\[
x_0^{r_x g r_y^{-1}}=x_0
\]
and thus $r_x gr_y^{-1}\in\dom\phi_{x_0}$. We have
\[\phi_{x_0}(r_x gr_y^{-1})=\left(r_x gr_y^{-1}\right)|_{x_0}.\]
But then
\[
(xw)^g=(xw)^{r_x^{-1}\cdot r_xgr_y^{-1}\cdot
r_y}=y\left(w^{h_x^{-1}\cdot\phi_{x_0}(r_xgr_y^{-1})\cdot
h_y}\right).
\]
We have used~(\ref{eq:restr3}) and~(\ref{eq:restr4}).

Thus we see that~(\ref{eq:nat2}) is true for the triple
$\left(\phi_{x_0}, T, C\right)$ and consequently, the action of
$G$ is defined by it. \qed\end{proof}

If the virtual endomorphism $\phi$ is induced by a partial
self-covering, then the set of actions defined by the triples
$(\phi, T, C)$ coincides with the set of standard actions, as
Propositions~\ref{pr:nat2} and~\ref{pr:natinv} show.

\begin{proposition}
\label{pr:nat2} Let $L=\{\ell_x\}_{x\in X}$ and $\LL:X\arr
f^{-1}(t)$ be a collection of paths and a bijection defining a
standard action. Take the virtual endomorphism $\phi_f$, induced
by $f$ and by a path $\ell$, connecting $t$ with $z_0\in
f^{-1}(t)$. For every $x\in X$, choose a path $\rho_x$ in $\M_1$,
starting at $z_0$ and ending at $\LL(x)$. Denote
$$r_x=f(\rho_x), \quad h_x=\ell\rho_x\ell_x^{-1}.$$ Then the
standard action is defined by the triple $\left(\phi_f,
\{r_x\}_{x\in X}, \{h_x\}_{x\in X}\right)$.
\end{proposition}

\begin{proof}
The set $\{r_x\}_{x\in X}$ is a right coset transversal for the
domain of  $\phi_f$, since $z_0^{r_x}=\LL(x)$ for every $x$. Let
$\gamma$ be an arbitrary loop at $t$. The path
$\gamma_x=f^{-1}(\gamma)\left[\LL(x)\right]$ ends in $\LL(y)$ for
$y=x^\gamma$. The element $r_xgr_y^{-1}$ belongs to $\dom\phi_f$,
since it is the loop $f(\rho_x\gamma_x\rho_y)$. The element
$\phi_f(r_xgr_y^{-1})$ is then the loop
$\ell\rho_x\gamma_x\rho_y\ell^{-1}$. Consequently
\[
h_x^{-1}\phi_f(r_xgr_y^{-1})h_y=\left(\ell\rho_x\ell_x^{-1}\right)^{-1}\ell\rho_x\gamma_x
\rho_y\ell^{-1}\left(\ell\rho_x\ell_x^{-1}\right)=\ell_x\gamma_x\ell_y^{-1},
\]
and Proposition~\ref{pr:nat} ends the proof. \qed\end{proof}


\begin{proposition}
\label{pr:natinv} Let $\phi$ be the virtual endomorphism of
$\pi(\M, t)$, defined by the partial self-covering $f$. Then for
every right coset transversal $T=\{r_x\}_{x\in X}$ and sequence
$C=\{h_x\}_{x\in X}$ there exists a collection of paths
$L=\{\ell_x\}_{x\in X}$ and a bijection $\LL:X\arr f^{-1}(t)$ such
that the respective standard action is the action defined by the
triple $(\phi, T, C)$.
\end{proposition}

\begin{proof}
Let $\ell$ be the path from $t$ to $z_0\in f^{-1}(t)$, inducing
together with $f$ the endomorphism $\phi$. For every $x\in X$ put
$\rho_x=f^{-1}(r_x)[z_0]$. Define $\LL(x)$ to be the end of
$\rho_x$. It follows from the fact that $r_x$ is a right coset
transversal of $\dom\phi$, that the defined map $\LL:X\arr
f^{-1}(t)$ is a bijection.

Take $\ell_x=h_x^{-1}\ell\rho_x$. Then
$h_x=\ell\rho_x\ell_x^{-1}$, $r_x=f(\rho_x)$ and
Proposition~\ref{pr:nat2} shows that the standard action defined
by the collection $L=\{\ell_x\}_{x\in X}$ and the bijection
$\LL:X\arr f^{-1}(t)$ coincides with the action defined by the
triple $(\phi, T, C)$. \qed\end{proof}

Finally, let us mention a description of the kernel of a
self-similar action (see~\cite{nek:stab} and~\cite{neksid}).

\begin{proposition}
\label{pr:kernel}
 The kernel of the action defined by a triple $(\phi, T, C)$ is equal to
$$
\mathcal{C}(\phi)=\bigcap_{n\ge 1}\bigcap_{g\in G}
g^{-1}\cdot\dom\phi^n\cdot g.
$$

In particular, the iterated monodromy group $\img{f}$ is the
quotient of the fundamental group $\pi_1(\M)$ by
$\mathcal{C}(\phi_f)$, where $\phi_f$ is the virtual endomorphism
induced by the self-covering.
\end{proposition}

\begin{proof}
Let $N$ be the kernel of the action defined by the triple. If
$g\in N$, then for every $x\in X$ we have
$r_xgr_x^{-1}\in\dom\phi$ and $h_x^{-1}\phi(r_xgr_x^{-1})h_x\in
N$. In particular, if we take $x=x_0$ such that
$r_{x_0}\in\dom\phi$, then we get that $g\in\dom\phi$, since
$r_{x_0}gr_{x_0}^{-1}\in\dom\phi$. Besides, we get that
$$h_{x_0}^{-1}\phi(r_{x_0}gr_{x_0}^{-1})h_{x_0}\in N,$$
hence
\[\phi(g)=(\phi(r_{x_0})^{-1}h_{x_0})\cdot\phi(g)\cdot(\phi(r_{x_0})^{-1}h_{x_0})^{-1}\in
N.\]

We have proved that $g\in N$ implies $g\in\dom\phi$ and
$\phi(g)\in N$. From this, by induction, we get that $\phi^n(g)\in
N$ and therefore, $\mathcal{C}(\phi)\ge N$.

On the other hand, if $g\in\mathcal{C}(\phi)$, then
$r_xgr_x^{-1}\in\dom\phi$ for every $x\in X$ and
$h_x^{-1}\phi(r_xgr_x^{-1})h_x\in\mathcal{C}(\phi)$. It follows by
induction on $n$ that the action of $g$ on $X^n$ is trivial. Hence
$\mathcal{C}(\phi)\le N$. \qed\end{proof}

\section{Contracting groups and expanding maps}
\label{s:contract}

\subsection{Contracting self-similar actions and their limit spaces}
\label{ss:contract}

We recall here some definitions and results of the
paper~\cite{nek:lim}.

\begin{defi}
\label{d:contr} A self-similar action of a group $G$ on the set
$\xs$ is said to be \emph{contracting} if there exists a finite
subset $\nuke\subset G$ such that for every $g\in G$ there exists
$n\in\N$ such that for every $v\in \xs$, $|v|>n$ the restriction
$g|_v$ belongs to $\nuke$.
\end{defi}

It is proved in~\cite{nek:lim}, that the property of an action to
be contracting depends only on the associated virtual endomorphism
$\phi$.

The minimal set $\nuke$ satisfying the conditions of
Definition~\ref{d:contr} is called the \emph{nucleus} of the
action. It follows from definition that if $h\in\nuke$ and $x\in
X$ then $h|_x$ belongs to $\nuke$, so we consider $\nuke$ as an
automaton.

Contraction of actions of finitely generated groups can be defined
using the contraction of the length of the group elements under
the action of the virtual endomorphism.

\begin{defi}
Let $\phi$ be a virtual endomorphism of a finitely generated group
$G$. Denote by $|g|$ the length of an element $g\in G$ with
respect to a fixed finite generating set $S=S^{-1}$, i.e., the
minimal length of a representation of $g$ as a product of the
elements of $S$.

Then the number
\[
\rho=\lim_{n\arr\infty}\sqrt[n]{\limsup_{g\in \dom\phi^n,
|g|\arr\infty} \frac{\left|\phi^n(g)\right|}{|g|}}
\]
is called the \emph{contraction coefficient} of the virtual
endomorphism $\phi$.
\end{defi}

The following proposition is proved in~\cite{nek:lim}.

\begin{proposition}
\label{pr:specr} The contraction coefficient of a virtual
endomorphism of a finitely generated group is finite and does not
depend on the choice of the generating set.

A level-transitive self-similar action of a finitely generated
group is contracting if and only if the contraction coefficient of
the associated virtual endomorphism is less than one.
\end{proposition}


Denote by $\xmo$ the set of all infinite to the left sequences of
the form $\ldots x_2x_1$ with the topology of the direct product
of discrete sets $\xmo=\cdots\times X\times X$.

\begin{defi}
\label{def:asympteq} Consider a self-similar action of a group $G$
on the set $\xs$. Two sequences $\ldots x_2x_1$ and $\ldots
y_2y_1$ are said to be \emph{asymptotically equivalent} with
respect to the action of $G$ if there exists a bounded sequence
$\{g_k\}_{k\geq 1}$ of elements of the group $G$, such that
$$
(x_kx_{k-1}\ldots x_1)^{g_k}=y_ky_{k-1}\ldots y_1
$$
for all $k\ge 1$.
\end{defi}

Here a sequence $\{g_k\}_{k\ge 1}$ is called \emph{bounded} if it
takes only a finite number of different values.

In other words, two sequences $\ldots x_2x_1$ and $\ldots y_2y_1$
are asymptotically equivalent if and only if the words $x_k\ldots
x_1$ and $y_k\ldots y_1$ stay on a uniformly bounded distance from
each other with respect to the action of $G$.

It is easy to see that the asymptotic equivalence is an
equivalence relation.


\begin{proposition}
\label{pr:limsp} Two sequences $\ldots x_2x_1$ and $\ldots y_2y_1$
are asymptotically equivalent with respect to a contracting action
if and only if there exists a directed path $\ldots e_2e_1$ in the
Moore diagram of the nucleus, for which the arrow $e_i$ is labeled
by $(x_i, y_i)$.
\end{proposition}

\begin{proof}
Let $\{g_k\}_{k\ge 1}$ be a bounded sequence of group elements
such that $$(x_k\ldots x_1)^{g_k}=y_k\ldots y_1.$$ There exists a
number $n_0$ such that $g_k|_{x_k\ldots x_{k-n}}$ belongs to the
nucleus for all $n\ge n_0$. Let $A_k$ be the set of all elements
of the nucleus of the form $g_{k+n}|_{x_{k+n}\ldots x_{k+1}}$.

It follows from the definitions that for every $a_k\in A_k$ we
have $$(x_k\ldots x_1)^{a_k}=y_k\ldots y_1$$ and $a_k|_{x_k}\in
A_{k-1}$.

All the sets $A_k$ are finite and non-empty. Thus, by a standard
argument, there exists a sequence $a_k\in A_k$ such that
$a_k|_{x_k}=a_{k-1}$ and $a_k(x_k)=y_k$. This sequence is the
necessary path in the Moore diagram. \qed\end{proof}

\begin{defi}
\label{def:limitspace} The quotient of the topological space
$\xmo$ by the asymptotic equivalence relation is called the
\emph{limit space} of the action, denoted $\lims{G}$. The
\emph{limit dynamical system} is the dynamical system $(\lims{G},
\si)$, where the map $\si:\lims{G}\arr\lims{G}$ is induced by the
shift
$$\sigma:\ldots x_2x_1\mapsto \ldots x_3x_2$$ on $\xmo$.
\end{defi}

The asymptotic equivalence is a shift-invariant, so that the map
$\si$ are well defined.

It is proved in~\cite{nek:lim} that the limit space $\lims{G}$ of
a contracting action is metrizable and finite dimensional. If the
group $G$ is finitely generated and level-transitive, then the
limit space is connected.

The following is also proved in~\cite{nek:lim}.

\begin{proposition}
\label{prop:limsp} The dynamical system $(\lims{G}, \si)$ is
uniquely determined, up to topological conjugacy, by the group $G$
and the conjugacy class of the virtual endomorphism $\phi$.
\end{proposition}

We have the following criterion for a metric space to be
homeomorphic to the limit space $\lims{G}$.

Another aspect of the limit space is that it can be represented as
a limit of the graphs of the action of $G$ on the levels $X^n$ of
the tree $\xs$. This can be formalized (see~\cite{nek:lim}), but
we will use this fact here only as an illustration (see
Figure~\ref{fig:jul1}).

\subsection{Iterated monodromy groups of expanding maps}
\label{ss:expand}

The space $\M$ in this subsection is a differentiable manifold and
the partial self-covering $f:\M_1\arr\M$ is smooth. If $\mu$ is a
Riemannian metric on a manifold, then by $d_\mu(x, y)$ we denote
the distance between the points $x$ and $y$ (i.e., the greatest
lower bound of the lengths of the piecewise smooth paths
connecting them).

\begin{defi}
\label{def:expanding} A map $f:\M_1\arr \M$, where $\M_1$ is an
open subset of $\M$,  is \emph{expanding} if there exist a
Riemannian metric $\mu$ on $\M$ and numbers $c>0, k>1$ such that
$\|Df^n \overrightarrow{v}\|_\mu\geq c\cdot
k^n\left\|\overrightarrow{v}\right\|_\mu$ for every $n\in\N$ and
every tangent vector $\overrightarrow{v}\in Tf^{-n}(\M)$.
\end{defi}

If the map $f$ is expanding and surjective, then it is a partial
self-covering of the space $\M$ and every piecewise smooth path of
length $l$ is mapped by $f^n$ onto a smooth path of length $\ge
c\cdot k^nl$.

\begin{defi}
\label{d:jul} The \emph{Julia set} of an expanding map $f$,
denoted $\mathcal{J}(f)$, is the set of the accumulation points of
the set $\cup_{n=0}^\infty f^{-n}(z_0)$, where $z_0\in\M$ is
arbitrary.

If $f$ is an expanding map, then we say that the set $\M$ is
\emph{complete on the Julia set} if every Cauchy subsequence of
the set $\cup_{n=0}^\infty f^{-n}(z_0)$ converges in $\M$.
\end{defi}

\begin{lemma}
If $f:\M_1\arr\M$ is an expanding map, then the Julia set
$\mathcal{J}(f)$ does not depend on the choice of the initial
point $z_0$ and $f\left(\mathcal{J}(f)\right)=\mathcal{J}(f)$,
$f^{-1}\left(\mathcal{J}(f)\right)=\mathcal{J}(f)$.
\end{lemma}

\begin{proof}
Let $z_1$ be another point in $\M$ and let $\gamma$ be a piecewise
smooth path, connecting $z_0$ with $z_1$. Then for every $x\in
f^{-n}(z_0)$ the end $y$ of the path $\gamma'=f^{-n}(\gamma)[x]$
belongs to $f^{-n}(z_1)$. The length of $\gamma'$  is not greater
than $c^{-1}k^{-n}(\mbox{{\it length}}(\gamma))$. Therefore
$d_\mu(x, y)\le c^{-1} k^{-n}d_\mu(z_0, z_1)$, what implies that
the set of accumulation points of $\cup_{n=0}^\infty f^{-n}(z_1)$
is equal to the set of the accumulation points of
$\cup_{n=0}^\infty f^{-n}(z_0)$, thus the Julia set does not
depend on the choice of $z_0$.

The set $\cup_{n=0}^\infty f^{-n}(z_0)$ is mapped by $f$ onto the
set $\cup_{n=0}^\infty f^{-n}(f(z_0))$, so that
$f\left(\cup_{n=0}^\infty f^{-n}(z_0)\right)=\cup_{n=0}^\infty
f^{-n}(z_0)\cup\{f(z_0)\}$, thus the Julia set is $f$-invariant.

We have also $f^{-1}\left(\cup_{n=0}^\infty
f^{-n}(z_0)\right)\cup\{z_0\}=\cup_{n=0}^\infty f^{-n}(z_0)$, so
that $f^{-1}\left(\mathcal{J}(f)\right)=\mathcal{J}(f)$.
\qed\end{proof}

\begin{defi}
We say that the fundamental group of a Riemannian manifold
$\mathcal{U}$ has \emph{finite balls} if for every two points $x,
y\in\mathcal{U}$ and every $R>0$ there exists only a finite number
of homotopy classes of paths of length $\le R$, connecting $x$
with $y$ in $\mathcal{U}$.
\end{defi}

The following lemma can be used to prove that the fundamental
group of a manifold has finite balls.

\begin{lemma}
\label{l:finhom} Let $\mathcal{U}$ be an arcwise connected
Riemannian manifold and let $\hat{\mathcal{U}}$ be its completion.
If the inclusion $\mathcal{U}\subset\hat{\mathcal{U}}$ induces an
isomorphism of the fundamental groups $\pi_1(\mathcal{U})$ and
$\pi_1\left(\hat{\mathcal{U}}\right)$, then the fundamental group
of $\mathcal{U}$ has finite balls.
\end{lemma}

\begin{proof}
The completion $\hat{\mathcal{U}}$ is a length space, i.e., the
distance between two its points is equal to the infimum of the
lengths of rectifiable paths, connecting them. Let
$\tilde{\mathcal{U}}$ be the universal covering of $\mathcal{U}$.
It is also a length space and its completion is the universal
covering $\widetilde{\mathcal{U}}$ of $\hat{\mathcal{U}}$. Then by
Hopf-Rinow Theorem (see~\cite{hopfrinow} and~\cite{bridhaefl}
p.~35), every bounded closed subset of $\widetilde{\mathcal{U}}$
is compact. Let $x_0$ be a preimage of the point $x\in\mathcal{U}$
in $\widetilde{\mathcal{U}}$. Every path of length not greater
than $R$ connecting $x$ with $y\in\mathcal{U}$ can be lifted to a
path of the same length, connecting $x_0$ with a preimage $y_i$ of
$y$. The homotopy class of the path is uniquely defined by the end
$y_i$ of the preimage. The set of all possible values of $y_i$ is
bounded and closed. Thus it is compact, i.e., finite.
\qed\end{proof}

\begin{theorem}
\label{th:limjul} Let $f:\M_1\arr \M$ be an expanding partial
self-covering map on $\M$. Suppose that  $\M$ is complete on the
Julia set $\mathcal{J}_{\M}(f)$ and that the Julia set has an open
arcwise connected neighborhood whose fundamental group has finite
balls. Then every standard action of $G=\img{f}$ on $\xs$ is
contracting and the limit dynamical system $\left(\lims{G},
\si\right)$ is topologically conjugated to the dynamical system
$\left(\mathcal{J}(f), f\right)$.
\end{theorem}

\begin{proof}
Let $L=\{\ell_x\}_{x\in X}$ be a set of paths defining a standard
action of $\img{f}$ on $\xs$. We may assume that the paths
$\ell_x$ are piecewise smooth. Let $l$ be the maximal length of
the paths from $L$.

The Julia set $\mathcal{J}(f)$ is equal to the set of the
accumulation points of $$\cup_{k=0}^\infty
f^{-k}(t)=T_t=\LL(\xs).$$

Let us define for all $v, u\in\xs$ a piecewise smooth path
$\ell(v; uv)$ starting at $\LL(v)$ and ending at $\LL(uv)$ by the
conditions:
\begin{enumerate}
\item $\ell(v; v)$ is the trivial path at $\LL(v)$;
\item $\ell(v; xv)=f^{-|v|}(\ell_x)\left[\LL(v)\right]$;
\item $\ell(w; uvw)=\ell(w; vw)\ell(vw; uvw)$.
\end{enumerate}

It easily follows from the definition of $\LL$ that the paths
$\ell(v; uv)$ are well defined. We have also that the length of
the path $\ell(v; uv)$ is not greater than
$c^{-1}k^{-|v|}\cdot(l+k^{-1}l+k^{-2}l+\cdots k^{-|u|+1}l)$.

Let $\ldots x_2x_1\in\xmo$. Each path from the sequence of the
paths
\[
\{\ell(x_n\ldots x_2x_1; x_m\ldots x_{n+1}x_n\ldots
x_2x_1)\}_{m\ge n},
\]
is a continuation of the previous one. In the limit we get a path
denoted $\ell(x_n\ldots x_2x_1; \ldots x_2x_1)$ of length not
greater than $c^{-1}l(k^{-n}+k^{-n-1}+\ldots)=
\frac{c^{-1}lk^{-n}}{1-k^{-1}}$.

Every $\ldots x_2x_1\in\xmo$ defines a sequence
$z_n=\LL(x_nx_{n-1}\ldots x_1)$, $n=0,1,\ldots$ of points of the
preimage tree $T_t$, i.e., a sequence of elements of the set
$\cup_{n=0}^\infty f^{-n}(t)$.

The length of the path $\ell(x_{n-1}\ldots x_1; x_nx_{n-1}\ldots
x_1)$, connecting $z_{n-1}$ with $z_n$ is not greater than
$c^{-1}k^{-n+1}l$. This implies that the sequence $z_n$ is a
Cauchy sequence in the metric space $(\M, d_\mu)$. Let us denote
its limit by $\LL(\ldots x_2x_1)$.  Then the path
$\ell(x_{n-1}\ldots x_2x_1; \ldots x_2x_1)$ starts at
$\LL(x_n\ldots x_1)$ and ends at $\LL(\ldots x_2x_1)$.

We have the following obvious properties of the map
$\LL:\xmo\arr\mathcal{J}(f)$:
\begin{equation}
\label{eq:LL1} d_\mu(\LL(v), \LL(\ldots x_2x_1v)) \leq
c^{-1}l\left(\sum_{i=0}^\infty k^{-|v|+i}\right)=
k^{-|v|}\frac{c^{-1}l}{1-k^{-1}} \end{equation} \begin{equation}
\label{eq:LL2} d_\mu(\LL(\ldots x_2x_1v), \LL(\ldots y_2y_1v))
\leq k^{-|v|}\frac{2lc^{-1}}{1-k^{-1}},
\end{equation}
for all $v\in\xs$ and $\ldots y_2y_1, \ldots x_2x_1\in\xmo$.

It follows from~(\ref{eq:LL2}) that the map
$\LL:\xmo\arr\mathcal{J}(f)$ is continuous.
Inequality~(\ref{eq:LL1}) implies that $\LL$ is onto.
Consequently, $\mathcal{J}(f)$ is compact as an image of the
compact space $\xmo$.

Since $f(\LL(x_nx_{n-1}\ldots x_1))=\LL(x_nx_{n-1}\ldots x_2)$ for
every $n$, we have
\begin{equation}
\label{eq:shift} f(\LL(\ldots x_2x_1))=\LL(\ldots x_3x_2).
\end{equation}

The spaces $\xmo$ and $\mathcal{J}(f)$ are compact, thus the map
$\LL$ is a quotient map (see~\cite{engel:outline} Theorem~9 on
p.~114). So it is sufficient to prove that the group $\img{f}$ is
contracting and that two points $\ldots x_2x_1, \ldots
y_2y_1\in\xmo$ are asymptotically equivalent if and only if their
$\LL$-images are equal.

Let us prove that the group $\img{f}$ is contracting. Let
$\mathcal{U}$ be an arcwise connected open neighborhood of
$\mathcal{J}(f)$ such that the fundamental group of $\mathcal{U}$
has finite balls. Consider an arbitrary infinite word $\ldots
x_2x_1\in\xmo$. The sequence $\{\LL(x_nx_{n-1}\ldots x_1)\}_{n=1,
2, \ldots}$ converges to the point $$\LL(\ldots
x_2x_1)\in\mathcal{J}(f)\subset\mathcal{U}.$$ It follows
from~(\ref{eq:LL1}) that there exists $n_0\in\N$, such that the
path $$\ell(x_{n_0}x_{n_0-1}\ldots x_1; wx_{n_0}x_{n_0-1}\ldots
x_1)$$ is inside the set $\mathcal{U}$  for every $w\in\xmo$.

In this way we cover the space $\xmo$ by cylindrical sets $\xmo
x_{n_0}\ldots x_1$, so that for any $w\in\xs x_{n_0}\ldots
x_1\cup\xmo x_{n_0}\ldots x_1$ and $v\in\xs x_{n_0}\ldots x_1$ the
path $\ell(v; w)$ belongs to $\mathcal{U}$. The space $\xmo$ is
compact, so we can choose a finite sub-cover $\{\xmo v_i\}_{i=1,
\ldots ,m}$, where $v_i\in\xs$. Denote the set $\{v_i\}_{i=1,
\ldots ,m}$ by $\mathsf{V}$. For a given $R>0$ define
$\mathsf{K}(R)$ to be the set of the elements of $\img{f}$ defined
by the loops of the form $\ell(\emp; v)\gamma\ell(\emp; u)^{-1}$,
where $v, u\in\mathsf{V}$ and $\gamma$ is a path in $\mathcal{U}$
of length not greater than $R$. Then the set $\mathsf{K}(R)$ is
finite for every $R$.

Let $\gamma$ be a loop at $t$, defining an element $g\in\img{f}$.
It follows from Proposition~\ref{pr:nat} that
$$g|_v=\ell(\emp;
v)\gamma_v\ell(\emp; u)^{-1},$$ where
$\gamma_v=f^{-|v|}(\gamma)\left[\LL(v)\right]$ and $u=v^\gamma$.
There exist, for $v$ long enough, words $v', u'\in\mathsf{V}$ such
that $v=w_1v'$ and $u=w_2u'$. Then
$$\ell(\emp; v)\gamma_v\ell(\emp; u)^{-1}=\ell(\emp; v')\ell(v'; w_1v')
\gamma_v\ell(u'; w_1u')^{-1}\ell(\emp; u')^{-1}.$$

The middle part $\ell(v'; w_1v')\gamma_v\ell(u'; w_1u')^{-1}$ of
the path is inside the set $\mathcal{U}$, if $v$ is long enough
(since then $\gamma_v$ is short), and its length is not greater
than
\begin{eqnarray*}
c^{-1}k^{-|v'|}(l+k^{-1}l+\cdots +k^{-|w_1|+1}l) &+& c^{-1}k^{-|v|}\mbox{{\it length}}(\gamma)+\\
c^{-1}k^{-|u'|}(l+k^{-1}l+\cdots +k^{-|w_2|+1}l) &<&
R_1+c^{-1}k^{-|v|}\mbox{{\it length}}(\gamma),
\end{eqnarray*}
where
$R_1=\max_{v'\in\mathsf{V}}\frac{k^{-|v'|}lc^{-1}}{1-k^{-1}}$.

So for all $v\in\xs$ long enough the restriction $g|_v$ belongs to
the set $\mathsf{K}(R_1+1)$, which is finite. Therefore, the
action of $\img{f}$ on $\xs$ is contracting.

Suppose that the points $\ldots x_2x_1, \ldots y_2y_1\in\xmo$ are
asymptotically equivalent. Then there exists a bounded sequence
$\{\gamma_n\}$ of loops at $t$ such that $(x_nx_{n-1}\ldots
x_1)^{\gamma_n}=y_ny_{n-1}\ldots y_1$ for every $n\geq 1$. Denote
by $m$ the maximal length of the paths $\gamma_n$. Let
$$
\gamma_n'=f^{-n}(\gamma_n)\left[\LL(x_nx_{n-1}\ldots x_1)\right].
$$
Then the end of $\gamma_n'$ is $\LL(y_ny_{n-1}\ldots y_1)$ and its
length will be not greater than $c^{-1}k^{-n}m$. Therefore
$$d_\mu(\LL(x_nx_{n-1}\ldots x_1), \LL(y_ny_{n-1}\ldots y_1))\le c^{-1}k^{-n}m,$$
so
\begin{eqnarray*}
\LL(\ldots x_2x_1) &=& \lim_{n\to\infty} \LL(x_nx_{n-1}\ldots
x_1)=\\
\lim_{n\to\infty}\LL(y_ny_{n-1}\ldots y_1) &=& \LL(\ldots y_2y_1).
\end{eqnarray*}

 Suppose now that $\LL(\ldots x_2x_1)=\LL(\ldots y_2y_1)$.
It follows from~(\ref{eq:shift}) that
$$
\LL(\ldots x_{n+2}x_{n+1})=\LL(\ldots y_{n+2}y_{n+1})
$$
for every $n$.

Then the path $$\gamma_n'=\ell(x_n\ldots x_1; \ldots
x_{n+2}x_{n+1}) \ell(y_n\ldots y_1; \ldots y_{n+2}y_{n+1})^{-1}$$
begins at $\LL(x_nx_{n-1}\ldots x_1)$ and ends in
$\LL(y_ny_{n-1}\ldots y_1)$. Its image under $f^n$ is the loop
$\gamma_n=\ell(\emp; \ldots x_2x_1)\ell(\emp; \ldots
y_2y_1)^{-1}$.

The path $\gamma_n$ is equal to
\begin{eqnarray*}
\ell(\emp; x_{n_1}x_{n_1-1}\ldots x_1) &&
\ell(x_{n_1}x_{n_1-1}\ldots x_1; \ldots x_2x_1)\cdot\\
\cdot\ell(y_{n_2}y_{n_2-1}\ldots y_1; \ldots y_2y_1)^{-1} &&
\ell(\emp; y_{n_2}y_{n_2-1}\ldots y_1)^{-1}
\end{eqnarray*}
for some $x_{n_1}x_{n_1-1}\ldots x_1,
y_{n_2}y_{n_2-1}\ldots y_1\in \mathsf{V}$. The middle part
$$\ell(x_{n_1}x_{n_1-1}\ldots x_1; \ldots x_2x_1) \ell(y_{n_2}y_{n_2-1}\ldots y_1; \ldots y_2y_1)^{-1}$$
is inside the set $\mathcal{U}$ and its length is not greater than
$\frac{c^{-1}l\left(k^{-n_1}+k^{-n_2}\right)}{1-k^{-1}}=R_2$, so
the path $\gamma_n$ belongs to a finite set $\mathsf{K}(R_2)$.
Then $(x_nx_{n-1}\ldots x_1)^{\gamma_n}=y_ny_{n-1}\ldots y_1$, and
the words $\ldots x_2x_1, \ldots y_2y_1$ are asymptotically
equivalent.

So the map $\LL:\xmo\arr\mathcal{J}_{\M}(f)$ induces a
homeomorphism of $\lims{G}$ with $\mathcal{J}(f)$. It follows
from~(\ref{eq:shift}) that this homeomorphism conjugates the shift
$\si$ with the map $f$.\qed\end{proof}

\section{Examples and applications}
\label{s:examps}

\subsection{Self-coverings and expanding endomorphisms of manifolds}

Here we consider the case when $f$ is defined on the whole space
$\M$. Then the respective virtual endomorphism $\phi_f$ of the
fundamental group $\pi_1(\M)$ is an isomorphism from a subgroup of
finite index $\dom\phi_f<\pi_1(\M)$ to $\pi_1(\M)$. Its inverse is
an injective endomorphism $f_{\#}:\pi_1(\M)\arr\pi_1(\M)$ induced
by the map $f$ (more pedantically, $f_{\#}$ is a map from
$\pi_1(\M, x)$ to $\pi_1(\M, t)$, where $x$ is a preimage of $t$,
but we identify $\pi_1(\M, x)$ with $\pi_1(\M, t)$ using a path
$\ell_x$ connecting $t$ to $x$).

The kernel of the iterated monodromy action of the fundamental
group $\pi_1(\M)$ is equal, by Proposition~\ref{pr:kernel}, to the
subgroup
$$N_f=\bigcap_{k\ge 1}\bigcap_{g\in G} g^{-1}\cdot f_{\#}^k\left(\pi_1\left(\M\right)\right)\cdot g.$$
The iterated monodromy group $\img{f}$ is isomorphic then to the
quotient $\pi_1(\M)/N_f$.

The following properties of expanding endomorphisms of Riemannian
manifolds where proved by M.~Shub and
J.~Franks~\cite{shub1,shub2}.

\begin{theorem}[M.~Shub, J.~Franks]
\label{th:shub} Suppose that the map $f:\M\arr\M$ on a compact
Riemannian manifold $\M$ is expanding. Then the following is true.
\begin{enumerate}
\item
 The map $f$ has a fixed point.
\item
 The universal covering space of $\M$ is diffeomorphic to $\R^n$.
\item
 The periodic points of $f$ are dense in $\M$.
\item
There exists a dense orbit of $f$ (i.e., the dynamical system
$(\M, f)$ is \emph{topologically transitive}).
\item
The fundamental group $\pi_1(\M)$ is a torsion free group of
polynomial growth.
\item
$$
\bigcap_{k\ge 1}f_{\#}^k\left(\pi_1(\M)\right)=\{1\}.
$$
\end{enumerate}
\end{theorem}

Theorem~\ref{th:limjul} and~\ref{th:shub} imply
\begin{theorem}
\label{th:expman} Let $f:\M\arr\M$ be an expanding map on a
compact manifold $\M$. Then the iterated monodromy group
$G=\img{f}$ is isomorphic to the fundamental group $\pi_1(\M)$.
Every its standard self-similar action on the tree $\xs$ is
contracting and the limit dynamical system $(\lims{G}, \si)$ is
topologically conjugate with the system $(\M, f)$.
\end{theorem}

\begin{proof}
The only thing we need to prove is that the set $\cup_{k\ge
1}f^{-1}(x)$ is dense in $\M$, i.e., that $\mathcal{J}(f)=\M$. But
this follows easily from the topological transitivity of $(\M,
f)$. Let $x_0\in\M$ be a point with a dense $f$-orbit. Then for
every $y\in\M$ and $\epsilon>0$ there exist $m$ and $n>m$ such
that $d(f^n(x_0), x)<\epsilon$ and $d(f^m(x_0), y)<\epsilon$, but
then
$$d(x', y)\le d(x', f^m(x))+d(f^m(x), y)<c^{-1}\cdot k^{m-n}\epsilon+\epsilon<(c^{-1}+1)\epsilon$$
for some $x'\in f^{-(n-m)}(x)$. Here $k>1$ and $c>0$ are the
constants from the definition of an expanding map. \qed\end{proof}

Theorem~\ref{th:expman} and Proposition~\ref{pr:natinv} imply the
following (see~\cite{shub1} Theorems~4 and~5).

\begin{theorem}[M.~Shub]
\label{th:shub2} The expanding map $f:\M\arr\M$ is uniquely
determined, up to topological conjugacy, by the action of $f_{\#}$
on its fundamental group $\pi_1(\M)$.
\end{theorem}

M.~Gromov, using his theorem on groups of polynomial growth, has
proved a conjecture of M.~Shub (see~\cite{shub2}
and~\cite{hirsch}), which describes all possible expanding
endomorphisms of Riemannian manifolds.

Let $L$ be a connected and simply connected nilpotent Lie group
and let $\mathop{\mathrm{Aff}}(L)$ be the group generated by the
left translations and automorphisms of the group $L$. Chose a
subgroup $G<\mathop{\mathrm{Aff}}(L)$ acting freely and discretely
on $L$. Suppose that the quotient $\M=L/G$ is compact. Then it is
a manifold. If an expanding endomorphism $\tilde f$ of the Lie
group $L$ conjugates $G$ to its subgroup, then $\tilde f$ induces
an expanding map $f:\M\arr\M$. Such map is called \emph{expanding
endomorphism} of the \emph{infranil-manifold} $\M$.

Note that an endomorphism of a Lie group is expanding if and only
if its derivative at $1$ is an expanding linear map.

\begin{theorem}[M.~Gromov]
Every expanding map of a compact manifold is topologically
conjugate to an expanding endomorphism of an infranil-manifold.
\end{theorem}

Explicit self-similar actions on the tree $\xs$ are interesting
from computational and dynamical points of view. They produce
faithful actions of the groups by finite-automatic automorphisms
of the tree $\xs$ and correspond to generalized numeration systems
on the group.

The space $\M$ is represented as a quotient of the space $\xmo$ by
the asymptotic equivalence relation, which is described in a
simple way by a finite graph (Proposition~\ref{pr:limsp}). This
gives an explicit symbolic finite presentation of the dynamical
system $(\M, f)$ (in the sense of M.~Gromov,
see~\cite{gro:hyperb,fried}). The images of the cylindrical sets
$\xmo v$ for $v\in\xs$ will define Markov partitions of $(\M, f)$
and define self-replicating tilings of the Lie groups.

Let us show this on several examples.

\paragraph{The adding machine.} Let $\mathbb{T}=\mathbb{R}/\mathbb{Z}$ be the circle.
The map $f:x\mapsto 2x$ induces a two-fold self-covering of the
circle $\mathbb{T}$. The fundamental group of the circle is
generated by the loop $\tau$ equal to the image of the segment
$[0, 1]$ in $\mathbb{R}/\mathbb{Z}$. Let us compute the standard
iterated monodromy action of the element $\tau$ on the tree $\{0,
1\}^*$.

Choose the base-point $t$ equal to $0$. It has two preimages:
itself and $1/2$. So we can take $\ell_0$ equal to the trivial
path in the point $0$ and $\ell_1$ equal to the image of the
segment $[0, 1/2]$.

The path $\tau$ has two preimages.  One is $\tau_0=\ell_1=[0,
1/2]$, another is $\tau_1=[1/2, 1]$.  Therefore $\tau$ acts on the
first level of the tree $\xs$ by the transposition. So
using~(\ref{eq:nat}) we get
\[
(0w)^\tau=1w, \quad (1w)^\tau=0w^\tau,
\]
since the path $\ell_0\tau_0 \ell_1^{-1}$ is trivial and
$\ell_1\tau_1\ell_0^{-1}$ is equal to $\tau$.

The recurrent definition of the transformation $\tau$ coincides
with the rules of adding 1 to a binary number. More precisely,
$(a_0\ldots a_n)^\tau=(b_0\ldots b_n)$ is equivalent to the
equality
$$
\left(a_0+a_1\cdot 2+\cdots +a_n\cdot 2^n\right)+1=
\left(b_0+b_1\cdot 2+\cdots+b_n\cdot 2^n\right)
$$
modulo $2^{n+1}$.

If we identify the infinite words $w=a_0a_1\ldots\in\{0,
1\}^\omega$ with the dyadic integers $\Phi(w)=a_0+a_1\cdot
2+a_2\cdot 2^2+\cdots$, then $\Phi(w^\tau)=\Phi(w)+1$.

The transformation $\tau$ is called the \emph{adding machine}. It
is an important example of a minimal dynamical system
(see~\cite{buescu}, for example).

\begin{figure}[ht]
\begin{center}
\includegraphics{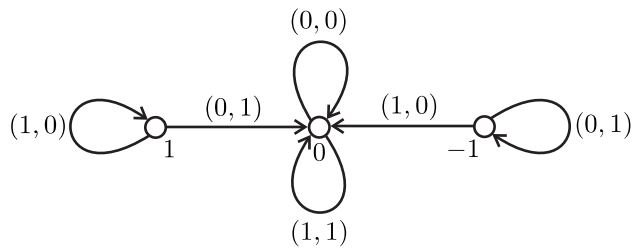}
\end{center}
\caption{The nucleus of the adding machine action} \label{fig:nuc}
\end{figure}

The map $f$ is obviously expanding and the described action of the
group $\Z$ on $\xs$ is contracting, since the respective virtual
endomorphism $\phi_f$ is the map $x\mapsto x/2$. The nucleus of
the action is the set $\{-1, 0, 1\}$. The Moore diagram of the
nucleus is shown on Figure~\ref{fig:nuc}.

The Moore diagram and Proposition~\ref{pr:limsp} show that two
sequences are asymptotically equivalent if and only if they are
equal or are of the form $\ldots 110v, \ldots 001v$, where
$v\in\xs$ is arbitrary, or of the form $\ldots 11, \ldots 00$. But
this is exactly the usual identification of the real binary
numbers, so the sequences $\ldots x_2x_1, \ldots y_2y_1\in\xmo$
are asymptotically equivalent if and only if
$$
x_1/2+x_2/4+\cdots+x_n/2^n+\cdots=y_1/2+y_2/4+\cdots+y_n/2^n+\cdots
(\mbox{mod\ } 1).
$$

It follows that the limit space is the circle $\R/\Z$ with the
shift map $\si(x)=2x (\mbox{mod\ } 1)$. We have returned back to
the original self-covering, what agrees with
Theorem~\ref{th:limjul}.

We see that the standard action of the iterated monodromy group on
the limit space defines an encoding of the space $\M=\R/\Z$ by
sequences of digits, which coincides with the binary numeration
system on $\R$.

\paragraph{The torus.}
The adding machine example can be generalized to $n$ dimensions.
Let $\mathbb{T}^n=\R^n/\Z^n$ be the $n$-dimensional torus. Let $A$
be an $n\times n$-matrix with integral entries and with
determinant equal to $d>1$.  Then the linear map $A:\R^n\arr\R^n$
induces an $d$-fold self-covering of the torus $\mathbb{T}^n$.

The fundamental group of the torus $\mathbb{T}^n$ is the free
Abelian group $\mathbb{Z}^n$. The iterated monodromy action is
faithful if and only if the group $N_A=\cap_{k\ge 1} A^k(\Z^n)$ is
trivial. This is equivalent to the condition that the no
eigenvalue of $A^{-1}$ is an algebraic integer (see~\cite{neksid}
Proposition~4.1 and~\cite{cuntz_rep} Proposition~10.1).

The corresponding standard iterated monodromy actions of $\Z^n$ on
rooted trees can be interpreted as numeration systems on $\Z^n$.
Consider a standard action of the group $\Z^n$ on $\xs$. Then it
follows from Proposition~\ref{pr:nat2} that there exists a coset
transversal $\{r_0, r_1, \ldots r_{d-1}\}$ for the subgroup
$A(\Z^n)$ such that
\[
(xw)^a=yw^{A^{-1}(a+r_x-r_y)},
\]
for all $x\in X$, $a\in\Z^n$ and $w\in\xs$, where $\{0, 1, 2,
\ldots, d-1\}=X$ and $y\in X$ is such that $a+r_x-r_y\in A(\Z^n)$.
We do not need the elements $h_x$, since $\phi_f$ in our case is
onto (see the remark after Definition~\ref{def:ssf}).

Consequently, $(x_0x_1\ldots x_m)^a=y_0y_1\ldots y_m$ is
equivalent to the condition that the elements
\[
\left(r_{x_0}+A\left(r_{x_1}\right)+A^2\left(r_{x_2}\right)+\cdots+A^m\left(r_{x_m}\right)\right)+a
\]
and
\[
r_{y_0}+A\left(r_{y_1}\right)+A^2\left(r_{y_2}\right)+\cdots+A^m\left(r_{y_m}\right)
\]
are equal modulo $A^{n+1}\left(\Z^n\right)$. So the self-similar
action of $\Z^n$ on the tree corresponds to an ``$A$-adic''
numeration system on $\Z^n$.

If the matrix $A$ is expanding, i.e., if all its eigenvalues have
absolute value greater than 1, then the series
\begin{equation}
\label{eq:znfractions}
\sum_{n=1}^\infty A^{-n}(r_{x_n})
\end{equation}
is convergent in $\R^n$ for every $\ldots x_2x_1\in\xmo$. In this
way the $A$-adic numeration system on $\Z^n$ extends to an
$A$-adic numeration system on $\R^n$.

It is proved in~\cite{nek:lim}, that in the case when $A$ is
expanding, the limit space of the iterated monodromy action is the
torus $\R^n/\Z^n$ and that the quotient map $\xmo\arr \R^n/\Z^n$
comes from this $A$-adic numeration system on $\R^n$ in a similar
way like for the case of the adding machine.

The set $\mathcal{T}$ of all possible sums of the
series~(\ref{eq:znfractions}) is called the \emph{digit tile}
defined by the $A$-adic numeration system. The $\Z^n$-translations
of the set $\mathcal{T}$ cover the space $\R^n$. In general the
translates can overlap, but often they form a tiling of the space
$\R^n$ (see a criterion in~\cite{nek:lim}).

The linear map $A$ maps every tile of such a tiling to a union of
$d$ tiles. The tilings with this properties are called
\emph{self-affine} (\emph{self-replicating tilings, rep-tilings,
digit-tilings}).

For every $k\ge 1$, the images of the sets
$A^{-k}(\mathcal{T}+r)$, $r\in\Z^n$ in the torus $\R^n/\Z^n$ form
a Markov partition of the dynamical system $(\R^n/\Z^n, A)$. They
are equal to the images of the cylindrical sets $\xmo v, v\in X^k$
with respect to the presentation of the limit space $\R^n/\Z^n$ as
a quotient of $\xmo$ by the asymptotic equivalence relation.

\begin{figure}
\begin{center}
  \includegraphics{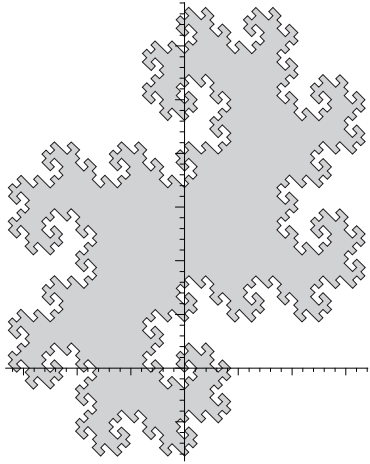}
\end{center}
\caption{The set of fractions bounded by the dragon curve.}
\label{fig:dragon}
\end{figure}

The set $\mathcal{T}$ often have a fractal boundary. One of most
known examples is the ``dragon curve'', shown on
Figure~\ref{fig:dragon}. It corresponds to the case
$A=\left(\begin{array}{cc} -1 & 1 \\ -1 & -1\end{array}\right)$
and the coset transversal $\{(0, 0), (1, 0)\}$. The respective
numeration system on $\R^2$ can be interpreted as a numeration
system on $\C$ with the base $(-1+i)$ and the digits $0, 1$. See
its discussion in~\cite{knuth}.

See the survey~\cite{vince:digtile} and bibliography in it for
properties of the digit tiles and their applications. The
question, when different $A$-adic expansions define the same point
in $\mathbb{R}^2$ was studied in~\cite{gilbert:three}.

\noindent
\paragraph{Heisenberg group.}
This example of a self-covering is from~\cite{shub1}. Let $L$ be
the group of lower triangular matrices
$
\left(
\begin{array}{ccc}
1 & 0 & 0\\
a & 1 & 0\\
c & b & 1
\end{array}
\right),
$
with $a, b, c\in\R$ and let $G$ be its subgroup of matrices with
$a, b, c\in\Z$. Then for all $p, q\in\Z$, the map
$$
f_{\#}:\left(
\begin{array}{ccc}
1 & 0 & 0\\
a & 1 & 0\\
c & b & 1
\end{array}
\right)\mapsto \left(
\begin{array}{ccc}
1 & 0 & 0\\
p\cdot a & 1 & 0\\
pq\cdot c & q\cdot b & 1
\end{array}
\right)
$$
is an automorphism of the group $L$, which maps $G$ to a subgroup
of index $p^2q^2$. The quotient $L/G$ is a three-dimensional
nil-manifold, and the map $f_{\#}$ induces its expanding
$p^2q^2$-fold self-covering.

Let us consider, for instance, the case $p=q=2$. Then by
Proposition~\ref{pr:natinv}, one of the standard self-similar
actions of the iterated monodromy group $G$ is defined by the pair
$(\phi, T)$, where $\phi$ is the virtual endomorphism
$\phi=f_{\#}^{-1}:\left(\begin{array}{ccc} 1 & 0 & 0\\
a & 1 & 0\\ c & b & 1 \end{array}\right)\mapsto \left(\begin{array}{ccc} 1 & 0 & 0\\
a/2 & 1 & 0\\ c/4 & b/2 & 1\end{array}\right)$ with the domain
$\left(\begin{array}{ccc} 1 & 0 & 0\\
2\Z & 1 & 0\\ 4\Z & 2\Z & 1 \end{array}\right)$ and $T$ is the
coset transversal
$\left\{\left(\begin{array}{ccc} 1 & 0 & 0\\
a & 1 & 0\\ c & b & 1 \end{array}\right): a, b\in\{0, 1\},
c\in\{0, 1, 2, 3\}\right\}$.

In the same way, as in the case of abelian groups, the images of
cylindrical sets $\xmo v$ are tiles of self-replicating tilings of
the Lie group $L$. See the paper~\cite{gotz:expand} for a
treatment of self-replicating tilings of groups. See also a
discussion of self-similar actions of nilpotent and solvable
groups in~\cite{neksid}.

\subsection{Rational functions on $\C$}
\subsubsection{Iterated monodromy group as a Galois group.}
\emph{Iterated monodromy group} of a rational function $f$ is by
definition the iterated monodromy group of the partial
self-covering $f:\M_1\arr \M$ of the the set
$\M=\hat{\C}\setminus\overline P$, where $\overline P$ is the
closure of the post-critical set of $f$ and $\M_1=f^{-1}(\M)$.

The following construction is due to R.~Pink (private
communication).

Let $f\in\C[z]$ be a polynomial over $\C$. For every $n\ge 1$
define a polynomial $F_n(z)=f^n(z)-t\in \C(t)[z]$ over the field
$\C(t)$ of rational functions, where $f^n(z)$ denotes the $n$th
iteration of $f$. Let $\Omega_n$ be the splitting field of $F_n$.
It is easy to see that $\Omega_n\subset\Omega_{n+1}$. It is a
classical fact, that the Galois group
$\mathop{Aut}(\Omega_n/\C(t))$ is isomorphic to the monodromy
group of the branched covering $f^n:\C\arr\C$ (see~\cite{forster}
Theorem~8.12), i.e., to the group of permutations of the set
$f^{-n}(z_0)$ induced by the action of the fundamental group
$\pi_1(\C\setminus P_n, z_0)$, where $P_n$ is the set of branching
points of $f^n$ and $z_0\notin P_n$ is an arbitrary point.

This implies the following interpretation of the iterated
monodromy group $\img{f}$.

\begin{proposition}
Let $f\in\C[z]$ be a post-critically finite polynomial. Then the
closure of the iterated monodromy group $\img{f}$ in the
automorphism group of the preimage tree is isomorphic to the
Galois group $\mathop{Aut}(\Omega/\C(t))$, where
$\Omega=\cup_{n\ge 1}\Omega_n$.
\end{proposition}

\subsubsection{Hyperbolic maps.}
The Julia set $\mathcal{J}(f)$ of a rational function $f\in\C(z)$
is the set of points $z\in\hat{\C}=\C\cup\{\infty\}$ such that the
set of functions $\{f^n: n\in \N\}$ is not normal on any
neighborhood of $z$ (see~\cite{milnor}). Since the set $\cup_{n\ge
0} f^{-n}(z)$ is dense in $\mathcal{J}(f)$ for every
$z\in\mathcal{J}(f)$ (see~\cite{milnor}, or~\cite{lyubich:top}),
our definition of the Julia set of an expanding map agrees with
the notion of the Julia set of a rational function.

\begin{defi}
A rational function $f\in\C(z)$ is said to be \emph{hyperbolic} if
it is expanding on a neighborhood of its Julia set.
\end{defi}

We have the following criterion, originally due to Fatou
(see~\cite{fatou}, see also a proof in~\cite{milnor},
Theorem~19.1).

\begin{theorem}
\label{th:hypfunc} A rational function $f\in\C(z)$ is hyperbolic
if and only if the closure of the postcritical set $\overline P$
does not intersect the Julia set $\mathcal{J}(f)$, or
equivalently, if and only if the orbit of every critical point
converges to an attracting cycle.
\end{theorem}

So, if a rational function is hyperbolic, then the post-critical
set has a finite number of accumulation points, which are all
outside the Julia set and the set
$\M=\hat{\C}\setminus\overline{P}$ is arcwise connected. Moreover,
the rational function is then expanding on the set $\M$ (see the
proof of Theorem~19.1 in~\cite{milnor}).

It is also easy to see that the Julia set has a neighborhood whose
fundamental group has finite balls. One can take, for instance,
the complement to the union of small closed disks around the
points of the attracting cycles (there exists only a finite number
of them) and around the post-critical points, which do not belong
to the already chosen disks.

Consequently, Theorem~\ref{th:limjul} implies

\begin{theorem}
\label{th:hypjul} Let $f\in\C(z)$ be a hyperbolic rational
function. Then every standard self-similar action of the iterated
monodromy group $\img{f}$ is contracting and the limit dynamical
system $\left(\lims{\img{f}}, \si\right)$ is topologically
conjugate to the dynamical system $(\mathcal{J}(f), f)$. In
particular, the limit space $\lims{\img{f}}$ is homeomorphic to
the Julia set of $f$.
\end{theorem}

Theorem~\ref{th:hypjul} provides a finite-to-one encoding of the
points of the Julia set $\mathcal{J}(f)$ by infinite sequences
over the alphabet $X$, which semi-conjugates the polynomial $f$ to
the Bernoulli shift. This and similar encodings where constructed
in~\cite{yacobson3,yacobson4,guckenheimer}, see
also~\cite{lyubich:top} p.~81--82.

Let us illustrate Theorem~\ref{th:hypjul} and computation of the
iterated monodromy groups (using Proposition~\ref{pr:nat}) by some
examples of polynomial mappings of degree $2$. All iterated
monodromy groups of these polynomials will act on the binary tree
$\xs=\{0, 1\}^*$. We will use here the permutational recursion,
described at the end of Subsection~\ref{ss:rooted}. Namely:
\begin{eqnarray*}
g=(g_0, g_1) &\Longleftrightarrow & (0w)^g=0w^{g_0} \mbox{ and } (1w)^g=1w^{g_1} \mbox{ for every } w\in\xs,\\
g=(g_0, g_1)\sigma &\Longleftrightarrow & (0w)^g=1w^{g_0} \mbox{
and } (1w)^g=0w^{g_1} \mbox{ for every } w\in\xs.
\end{eqnarray*}
Here $\sigma$ denotes the ``switch'' $(0w)^\sigma=1w,
(1w)^\sigma=0w$. For example, the adding machine transformation
$\tau$ is defined in this notation by the recurrent formula
$\tau=(1, \tau)\sigma$.

\paragraph{The adding machine as  $\img{z^2}$. } The polynomial $z^2$ defines on the circle
$\{z\in\C: |z|=1\}$ a self-covering, conjugate to the
self-covering $x\mapsto 2x$ of the circle $\R/\Z$. The
post-critical set of $z^2$ is $\{0, \infty\}$, so the fundamental
group of the space $\M=\hat{\C}\setminus P$ is generated by a loop
around the circle, and the computation of the standard action of
the iterated monodromy group of $z^2$ repeats the computation of
the iterated monodromy group of the self-covering $x\mapsto 2x$ of
the circle $\R/\Z$.

The polynomial $z^2$ is clearly hyperbolic (both critical points
$0, \infty$ are attracting fixed points). As we have already seen,
the limit space of the adding machine action is the circle, what
agrees with the fact that the Julia set of the polynomial $z^2$ is
the circle $\{z : |z|=1\}$.

\paragraph{Computation of $\img{z^2-1}$.}
The critical points of the polynomial $z^2-1$ are $\infty$ and
$0$. The infinity is a fixed point, and the orbit of $0$ is
$0\mapsto -1\mapsto 0$, so the post-critical set is $P=\{0, -1,
\infty\}$. The cycle $\{0, -1\}$ is attracting, since $0$ is a
critical point (so it is even \emph{super-attracting}).

\begin{figure}[h]
\begin{center}
\includegraphics{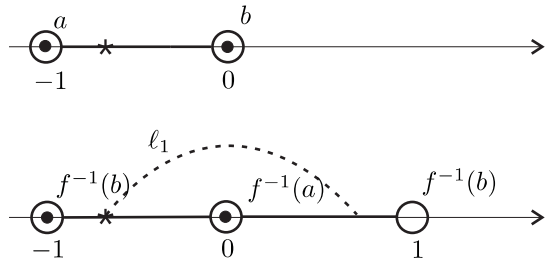}
\caption{Computation of the group $\img{z^2-1}$}\label{fig:min}
\end{center}
\end{figure}

Choose a basepoint $t=\frac{1-\sqrt{5}}{2}$ (denoted by a star on
Figure~\ref{fig:min}). It has two preimages: itself, and $-t$. Let
$\ell_0$ be the trivial path at $t$ and let $\ell_1$ be the path,
connecting $t$ and $-t$ as on the lower part of
Figure~\ref{fig:min}. Let $a$ and $b$ be the elements of
$\img{z^2-1}$, defined by the loops in positive direction around
$-1$ and $0$ respectively, shown on the upper part of the figure.

The preimages of the loops $a$ and $b$ are shown on the lower part
of Figure~\ref{fig:min}. We have
$$
a=(b, 1)\sigma, \quad b=(a, 1).
$$

\begin{figure}[h]
\begin{center}
\includegraphics{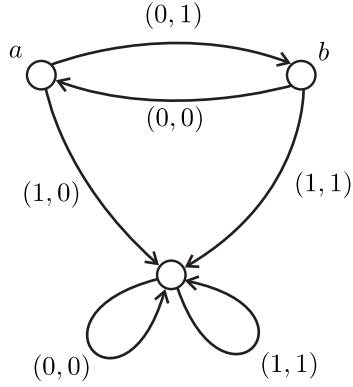}
\end{center}
\caption{The automaton generating the group
$\img{z^2-1}$.}\label{fig:moor1}
\end{figure}

Thus, the group $\img{z^2-1}$ is generated by the automaton with
the Moore diagram shown on Figure~\ref{fig:moor1}.

In the papers~\cite{zukgrigorchuk:3st,zukgrigorchuk:3stsp} the
following properties of the group $\img{z^2-1}$ are proved.

\begin{theorem}[R.~Grigorchuk, A.~\.Zuk]
\label{th:propimgmin}
 The group $\img{z^2-1}$
\begin{enumerate}
\item is torsion free;
\item has exponential growth (actually, the semigroup generated by $a$ and $b$ is free);
\item is just non-solvable, i.e., every its proper quotient is solvable;
\item
has solvable word and conjugacy problems;
\item
has no free non-abelian subgroups of rank 2.
\end{enumerate}
\end{theorem}

Theorem~\ref{th:propimgmin} is proved using the methods developed
during the study of groups acting on rooted trees, in particular
branch groups (see~\cite{grigorchuk:branch}) and of the just
non-solvable group of A.~Brunner, S.~Sidki and A.~Vieira
in~\cite{bsv:jns}. The properties of the group $\img{z^2-1}$ are
very similar to the properties of the group from~\cite{bsv:jns},
which is also a subgroup of the pro-finite completion of
$\img{z^2-1}$.

The following theorem is a result of L.~Bartholdi.

\begin{theorem}
\label{th:presimgmin}
 The group $\img{z^2-1}$ has the following presentation by defining relations:
$$
\img{z^2-1}=\left\langle a, b\left| \left[\left[a^{2^k},
b^{2^k}\right], b^{2^k}\right], \left[\left[b^{2^k},
a^{2^{k+1}}\right], a^{2^{k+1}}\right], k\ge
0\right.\right\rangle.
$$
\end{theorem}
Here $[x, y]=x^{-1}y^{-1}xy$.

The presentation of the group $\img{z^2-1}$ is similar to the
presentation of the Grigorchuk group, which was constructed by
I.~Lysionok in~\cite{lysionok:pres} and generalized for many other
contracting groups by L.~Bartholdi~\cite{bartholdi:lpres}.

Another interesting property of the group $\img{z^2-1}$ is its
amenability, proved by B.~Vir\'ag and L.~Bartholdi
(see~\cite{barthvirag}), using self-similar random walks. It was
proved before in~\cite{zukgrigorchuk:3stsp} that the group
$\img{z^2-1}$ does not belong to the class of sub-exponentially
amenable groups, i.e., can not be constructed from groups of
sub-exponential growth using the group-theoretic operations,
preserving amenability (passing to subgroups, quotients,
extensions and direct limits). The group $\img{z^2-1}$ is the
first example of an amenable group of this sort.

\begin{figure}[ht]
  \begin{center}
    \includegraphics{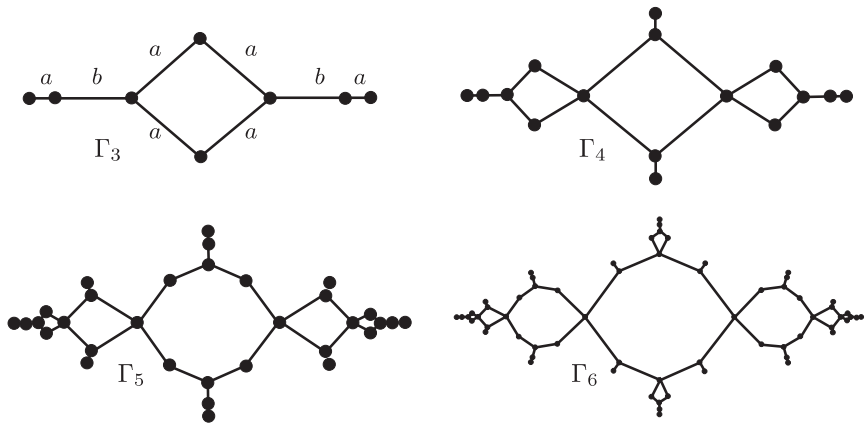}
\end{center}
  \caption{The Schreier graphs $\G_n(\img{z^2-1}, \{a, b\})$ for $3\le n\le 6$.}
  \label{fig:g2shr}
\end{figure}

\begin{figure}[ht]
    \begin{center}
      \includegraphics{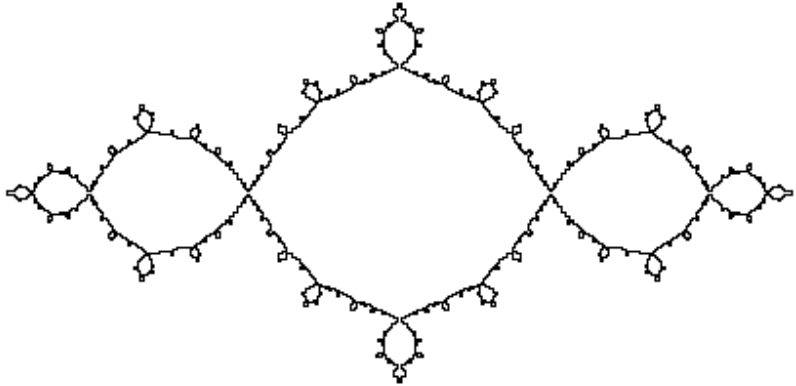}
    \end{center}
    \caption{The Julia set of the polynomial $z^2-1$.}
    \label{fig:jul1}
  \end{figure}

Figure~\ref{fig:g2shr} shows the Schreier graphs of action of the
group $\img{z^2-1}$ on the levels of the tree $\xs$. The shape of
these Schreier graphs was described by L.~Bartholdi
(see~\cite{bgn}). The Julia set of the polynomial $z^2-1$ is shown
on Figure~\ref{fig:jul1}.

\subsubsection{Sub-hyperbolic maps}
\label{sss:shm}

The set of hyperbolic functions is a subset of a more general
class of \emph{sub-hyperbolic} functions. A function is
sub-hyperbolic (see~\cite{milnor}) if it is expanding respectively
to some \emph{orbifold} metric on a neighborhood of the Julia set.
An analog of Theorem~\ref{th:hypfunc} is the following criterion
(see~\cite{milnor}).

\begin{theorem}
A rational function is sub-hyperbolic if and only if every orbit
of a critical point is either finite, or converges to an
attracting cycle.
\end{theorem}
In particular, every post-critically finite function (i.e., a
function for which the post-critical set $P$ is finite) is
sub-hyperbolic.

Theorem~\ref{th:hypjul} holds also for the sub-hyperbolic
functions. The proof of Theorem~\ref{th:limjul} (with small
modifications) is also valid for the case of orbifold coverings,
defined by sub-hyperbolic functions. More details will appear in a
subsequent paper, where iterated monodromy groups of orbifold
coverings will be studied.

Let us consider some examples

\paragraph{The dihedral group as $\img{z^2-2}$. }
The orbit of the finite critical point is $0\mapsto -2\mapsto
2\mapsto 2$. The post-critical set is $\{-2, 2, \infty\}$. Take,
for instance, the base point $t=0$ and connect it with the
preimages $\pm\sqrt{2}$ by straight segment. The fundamental group
of the space $\M=\C\setminus\{-2, 2\}$ is generated by a small
loop $a$ around $-2$, which is connected to $t$ by a straight
segment, and by a small loop $b$ around $2$, which is connected to
$t$ by a straight segment (both loops go around the points in the
positive direction). Computation of the standard action shows that
the respective generators of the iterated monodromy group
$\img{z^2-2}$ are defined by the recursion
$$
a=(1, 1)\sigma,\quad b=(b, a).
$$

It follows from the formula that $a^2=1$ and $b^2=(b^2, a^2)=(b^2,
1)=1$. So, the elements $a$ and $b$ are both of order 2. They
generate the infinite dihedral group $\mathbb{D}_\infty$
(see~\cite{grineksu_en} and~\cite{bgn}). The Schreier graph of
this group on the level $X^n$ is a chain of edges of length
$2^n-1$. It follows that the limit space is homeomorphic to the
real segment. This agrees with the fact that the Julia set of the
polynomial $z^2-2$ is the segment $[-2, 2]$.

This example can be generalized to the Chebyshev polynomials
$T_n(z)=\cos(n\arccos z)$, which all have the iterated monodromy
group isomorphic to $\mathbb{D}_\infty$ and Julia set $[-1, 1]$
(see~\cite{bgn}).

The Chebyshev polynomials are the only polynomials with the Julia
set $[-1, 1]$ (see~\cite{beardon:iter} Theorem~1.4.1).

\paragraph{The sphere and example of Latt\`es.}
If every critical point $z_0$ of a rational function $f$ is
pre-periodic (i.e., if $f^m(z_0)=f^k(z_0)$ for some $m\ne k$, but
$f^n(z_0)\ne z_0$ for all $n\ge 1$), then the Julia set of $f$ is
the whole sphere $\hat{\C}$ (see, for example Theorem~1.24
in~\cite{lyubich:top} or Theorem~9.4.4 in~\cite{beardon:iter}).

This is the case, for example, for the rational functions of
S.~Latt\`es~\cite{lattes:iter}. Let $\Gamma$ be a lattice in $\C$,
and let $\alpha$ be such that $\alpha\Gamma\subset\Gamma$. The
Weierstrass elliptic function $\wp$ for the lattice $\Gamma$ is
defined as
$$
\wp(z)=\frac{1}{z^2}+\sum_{\omega\in\Gamma\setminus\{0\}}\left[\frac{1}{(z+\omega)^2}-\frac{1}{\omega^2}\right].
$$
It induces a two-fold branched covering of the sphere $\hat{\C}$
by the torus $\C/\Gamma$. The map $z\mapsto \alpha z$ defines an
$|\alpha|^2$-covering of the torus $\C/\Gamma$ and respectively,
an $|\alpha|^2$-fold branched covering of the sphere $\hat{\C}$,
defined by a rational function $f$ of degree $|\alpha|^2$ such
that
$$
\wp(\alpha z)=f(\wp(z)).
$$

 For example, for $\alpha=2$ the function $f$ is
$$
f(z)=\frac{z^4+\frac{g_2}{2}z^2+2g_3z+\frac{g_2^2}{16}}{4z^3-g_2z-g_3},
$$
(see~\cite{beardon:iter} p.~74), where $g_2=60s_4$ and $g_3=140
s_6$ for $$s_m=\sum_{\omega\in\Gamma, \omega\ne 0} \omega^{-m}.$$

A pair $(g_2, g_3)$ is realized by a lattice $\Gamma$ if and only
if $g_2^3-27g_3^2\ne 0$ (see~\cite{lang:ellipt}, p.~39). In
particular, there exists a lattice $\Gamma$ such that $g_3=0$ and
$g_2=4$, so that
\begin{equation}
\label{eq:lattes} f(z)=\frac{(z^2+1)^2}{4z(z^2-1)}.
\end{equation}

For the case of the lattice $\Gamma=\Z[i]$ we have $g_3=0$, thus
$f(z)=\frac{(z^2+g_2/4)^2}{4z(z^2-g_2/4)}$, which is also
conjugate to~(\ref{eq:lattes}) (the conjugating map is
$t(z)=\frac{2z}{\sqrt{g_2}}$).

\begin{proposition}
\label{pr:lattes} Let $\Gamma$ be a lattice in $\C$ and let
$\alpha\in\C$ be such that $\alpha\Gamma\subset \Gamma$ and
$|\alpha|\ne 1$. Let a rational function $f\in\C(z)$ be such that
$$\wp(\alpha z)=f(\wp(z)).$$ Then the iterated monodromy group
$\img{f}$ is isomorphic to the group of affine transformations
$(-1)^kz+\omega$, where $k\in\Z, \omega\in\Gamma$. The associated
virtual endomorphism is the map $(-1)^kz+\omega\mapsto
(-1)^kz+\alpha^{-1}\omega$.
\end{proposition}

\begin{proof}
The Weierstrass function $\wp$ defines a two-fold branched
covering $\wp_0:\C/\Gamma\arr\hat{\C}$. The transformation
$\alpha:z\mapsto \alpha z$ defines a $|\alpha|^2$-fold
self-covering $f_0$ of the torus $\C/\Gamma$. The iterated
monodromy group $\img{f_0}$ is the group $\Gamma$ with the virtual
endomorphism $\phi_0:\omega\mapsto\alpha^{-1}\omega$.

Take the basepoint $t$ in $\hat{\C}$ not equal to a ramification
point of the function $\wp_0$. Then the point $t$ has two
preimages $z_0$ and $-z_0\in\C/\Gamma$ under $\wp_0$. The $n$th
level of the preimage tree $T_{z_0}$ of the point $z_0$ with
respect to the map $f_0$ is the image of $\alpha^{-n}(\Gamma+z_0)$
in the quotient $\C/\Gamma$. The image of
$\alpha^{-n}(\Gamma-z_0)$ in $\C/\Gamma$ is the $n$th level of the
preimage tree of the point $-z_0$.  Since
$f^n(\wp(z))=\wp(\alpha^n z)$, the $n$th level of the preimage
tree $T_t$ of the point $t$ respectively to the map $f$ is
$\wp_0\left(\alpha^{-n}\left(\Gamma+z_0\right)\right)=\wp_0\left(\alpha^{-n}\left(\Gamma-z_0\right)\right)$.

Let $\gamma$ be a loop based at $t$. Its preimages under $\wp_0$
are either a loop $\gamma_0$ based at $z_0$ and the symmetrical
loop $-\gamma_0$ based at $-z_0$ or a path $\gamma_0$ from $z_0$
to $-z_0$ and the symmetrical path $-\gamma_0$ from $-z_0$ to
$z_0$. In the first case every lift of the path $\gamma_0$ to the
universal cover $\C$ of the torus $\C/\Gamma$ connects the point
$z_0+\omega$, $\omega\in\Gamma$ to the point $z_0+\omega+a$ for
some fixed $a\in\Gamma$. Similarly, $-\gamma_0$ will connect the
point $-z_0-\omega-a$ to the point $-z_0-\omega$. In the second
case every lift of the path $\gamma_0$ to the universal covering
$\C$ of $\C/\Gamma$ connects the point $z_0+\omega$,
$\omega\in\Gamma$ to the point $-z_0+\omega+a$ for some fixed
$a\in\Gamma$ and every lift of the path $-\gamma_0$ connects the
point $-z_0-\omega-a$ with the point $-z_0-\omega$.

Consequently, in the first case $\gamma$ acts on the $n$th level
of the preimage tree $T_t$ by the transformation
$$
\wp_0\left(\alpha^{-n}(z_0+\omega)\right)\mapsto\wp_0\left(\alpha^{-n}(z_0+\omega+a)\right)
$$
and in the second case by the transformation
$$
\wp_0\left(\alpha^{-n}(z_0+\omega)\right)\mapsto\wp_0\left(\alpha^{-n}(-z_0+\omega+a)\right)=
\wp_0\left(\alpha^{-n}(z_0-\omega-a)\right).
$$
This implies immediately the statement of the proposition.
\qed\end{proof}

\bibliographystyle{plain}
\bibliography{mymath,nekrash}

\def\cprime{$'$}
\begin{thebibliography}{10}

\bibitem{shub2}
M.~Shub.
\newblock Expanding maps.
\newblock In {\em Global Analysis}, volume~14 of {\em Proc. Sympos. Pure
  Math.}, pages 273--276. American Math. Soc., Providence, Rhode Island, 1970.

\bibitem{gro:hyperb}
Mikhael Gromov.
\newblock Hyperbolic groups.
\newblock In S.~M. Gersten, editor, {\em Essays in Group Theory}, number~8 in
  M.S.R.I. Pub., pages 75--263. Springer, 1987.

\bibitem{forster}
Otto Forster.
\newblock {\em Lectures on {Riemann} surfaces}, volume~81 of {\em Graduate
  Texts in Mathematics}.
\newblock New York -- Heidelberg -- Berlin: Springer-Verlag, 1981.

\bibitem{vince:digtile}
Andrew Vince.
\newblock Digit tiling of {E}uclidean space.
\newblock In {\em Directions in mathematical quasicrystals}, pages 329--370.
  Amer. Math. Soc., Providence, RI, 2000.

\bibitem{lang:ellipt}
Serge Lang.
\newblock {\em {Elliptic functions. Second edition.}}, volume 112 of {\em
  Graduate Texts in Mathematics}.
\newblock {Springer-Verlag, New York etc.}, 1987.

\bibitem{engel:outline}
Ryszard Engelking.
\newblock {\em Outline of general topology}.
\newblock Amsterdam: North-Holland Publishing Company, 1968.

\bibitem{fatou}
Pierre Fatou.
\newblock Sur les \'equations fonctionnelles.
\newblock {\em Bull. Soc. Math. France}, 48:33--94, 1920.

\bibitem{barthvirag}
L.~Bartholdi and B.~Vir\'ag.
\newblock Amenability via random walks.
\newblock submitted.

\bibitem{eil}
Samuel Eilenberg.
\newblock {\em Automata, Languages and machines}, volume~{A}.
\newblock Academic Press, New York, London, 1974.

\bibitem{grigorchuk:branch}
Rostislav~I. Grigorchuk.
\newblock Just infinite branch groups.
\newblock In Aner Shalev, Marcus~P.~F. {du Sautoy}, and Dan Segal, editors,
  {\em New horizons in pro-$p$ groups}, volume 184 of {\em Progress in
  Mathematics}, pages 121--179. Birkh{\"a}user Verlag, Basel, etc., 2000.

\bibitem{sidki_monogr}
Said~N. Sidki.
\newblock {\em Regular Trees and their Automorphisms}, volume~56 of {\em
  Monografias de Matematica}.
\newblock IMPA, Rio de Janeiro, 1998.

\bibitem{fried}
David~L. Fried.
\newblock Finitely presented dynamical systems.
\newblock {\em Ergod. Th. Dynam. Sys.}, 7:489--507, 1987.

\bibitem{guckenheimer}
J.~Guckenheimer.
\newblock Endomorphisms of the {Riemann} sphere.
\newblock In {\em Global Analysis}, volume~14 of {\em Proc. Sympos. Pure
  Math.}, pages 95--123. American Math. Soc., Providence, Rhode Island, 1970.

\bibitem{massey:algtop}
William~S. Massey.
\newblock {\em A basic course in algebraic topology}, volume 127 of {\em
  Graduate Texts in Mathematics}.
\newblock Springer-Verlag, New York etc., 1991.

\bibitem{hirsch}
M.~W. Hirsch.
\newblock Expanding maps and transformation groups.
\newblock In {\em Global Analysis}, volume~14 of {\em Proc. Sympos. Pure
  Math.}, pages 125--131. American Math. Soc., Providence, Rhode Island, 1970.

\bibitem{sush:avt_en}
Vitali{\u\i}~I. Sushchansky.
\newblock Groups of automatic permutations.
\newblock {\em Dop. NAN Ukrainy}, (6):47--51, 1998.
\newblock (in Ukrainian).

\bibitem{gilbert:three}
W.~J. Gilbert.
\newblock Complex numbers with three radix expansions.
\newblock {\em Can.~J. Math.}, 34(6):1335--1348, 1982.

\bibitem{nek:iterat}
Volodymyr~V. Nekrashevych.
\newblock Virtual endomorphisms of groups.
\newblock {\em Algebra and Discrete Mathematics}, 1(1):96--136, 2002.

\bibitem{shub1}
M.~Shub.
\newblock Endomorphisms of compact differentiable manifolds.
\newblock {\em Am. J. Math.}, 91:175--199, 1969.

\bibitem{grigorchuk:milnor_en}
Rostislav~I. Grigorchuk.
\newblock On the {M}ilnor problem of group growth.
\newblock {\em Dokl.\ Akad.\ Nauk SSSR}, 271(1):30--33, 1983.

\bibitem{bartholdi:lpres}
Laurent Bartholdi.
\newblock Endomorphic presentations of branch groups.
\newblock {\em Journal of Algebra}, 268(2):419--443, 2003.

\bibitem{grigorchuk:80_en}
Rostislav~I. Grigorchuk.
\newblock On {Burnside's} problem on periodic groups.
\newblock {\em Funtional Anal.\ Appl.}, 14(1):41--43, 1980.

\bibitem{lyubichminsk}
Mikhail Lyubich and Yair Minsky.
\newblock Laminations in holomorphic dynamics.
\newblock {\em J. Differ. Geom.}, 47(1):17--94, 1997.

\bibitem{bridhaefl}
M.~R. Bridson and A.~Haefliger.
\newblock {\em Metric spaces of non-positive curvature}, volume 319 of {\em
  Grundlehren der Mathematischen Wissenschaften}.
\newblock Springer, Berlin, 1999.

\bibitem{buescu}
Jorge Buescu and Ian Stewart.
\newblock Liapunov stability and adding machines.
\newblock {\em Ergodic Theory Dynamical Systems}, 15:1--20, 1995.

\bibitem{nek:lim}
Volodymyr~V. Nekrashevych.
\newblock Limit spaces of self-similar group actions.
\newblock preprint, Geneva University, available at
  \texttt{http://www.unige.ch/math/biblio/preprint/2002/limit.ps}, 2002.

\bibitem{lyubich:top}
M.Yu. Lyubich.
\newblock The dynamics of rational transforms: the topological picture.
\newblock {\em Russ. Math. Surv.}, 41(4):43--117, 1987.

\bibitem{milnor}
John~W. Milnor.
\newblock {\em Dynamics in one complex variable. {I}ntroductory lectures.}
\newblock Wiesbaden: Vieweg, 1999.

\bibitem{gotz:expand}
G.~Gelbrich.
\newblock Self-similar tilings and expanding homomorphisms of groups.
\newblock {\em Arch. Math.}, 65(6):481--491, 1995.

\bibitem{lysionok:pres}
Igor~G. Lysionok.
\newblock A system of defining relations for the {Grigorchuk} group.
\newblock {\em Mat.\ Zametki}, 38:503--511, 1985.

\bibitem{grigorchuk:semigr_en}
Rostislav~I. Grigorchuk.
\newblock On cancellative semigroups of polynomial growth.
\newblock {\em Mat. Zametki}, 43(3):305--319, 1988.

\bibitem{knuth}
Donald~E. Knuth.
\newblock {\em The art of computer programming, {V}ol 2, {S}eminumerical
  Algorithms}.
\newblock Addison-Wesley Publishing company, 1969.

\bibitem{gupta-sidki_group}
Narain~D. Gupta and Said~N. Sidki.
\newblock On the {Burnside} problem for periodic groups.
\newblock {\em Math.~Z.}, 182:385--388, 1983.

\bibitem{yacobson4}
M.~V. Yacobson.
\newblock Markov partitions for rational endomorphisms of the {Riemann} sphere.
\newblock In {\em Multicomponent random systems}, pages 381--396. Dekker, New
  York, 1980.

\bibitem{dh:plike}
Adrien Douady and John~Hamal Hubbard.
\newblock On the dynamics of polynomial-like mappings.
\newblock {\em {Ann. Sci. \'Ec. Norm. Sup\'er. IV. S\'er.}}, 18:287--343, 1985.

\bibitem{hopfrinow}
H.~Hopf and W.~Rinow.
\newblock {\"Uber den Begriff der vollst\"andigen differentialgeometrischen
  Fl\"ache}.
\newblock {\em Comment. Math. Helv}, 3:209--225, 1932.

\bibitem{bgn}
L.~Bartholdi, R.~Grigorchuk, and V.~Nekrashevych.
\newblock From fractal groups to fractal sets.
\newblock In Peter Grabner and Wolfgang Woess, editors, {\em {Fractals in Graz
  2001. Analysis -- Dynamics -- Geometry -- Stochastics}}, pages 25--118.
  {Birkh\"auser Verlag. Basel. Boston. Berlin.}, 2003.

\bibitem{beardon:iter}
Alan~F. Beardon.
\newblock {\em {Iteration of rational functions. Complex analytic dynamical
  systems}}, volume 132 of {\em Graduate Texts in Mathematics}.
\newblock {Springer-Verlag. New York etc.}, 1991.

\bibitem{lattes:iter}
S.~Latt{\`e}s.
\newblock Sur l'it\'eration des substitutions rationelles et les fonctions de
  {P}oincar\'e.
\newblock {\em C.~R. Acad. Sci. Paris}, 166:26--28, 1918.

\bibitem{nek:stab}
Volodymyr~V. Nekrashevych.
\newblock Stabilizers of transitive actions on locally finite graphs.
\newblock {\em Int.~J. of Algebra and Computation}, 10(5):591--602, 2000.

\bibitem{grineksu_en}
Rostislav~I. Grigorchuk, Volodymyr~V. Nekrashevich, and Vitali{\u\i}~I.
  Sushchanskii.
\newblock Automata, dynamical systems and groups.
\newblock {\em Proceedings of the Steklov Institute of Mathematics},
  231:128--203, 2000.

\bibitem{bsv:jns}
Andrew~M. Brunner, Said~N. Sidki, and Ana.~C. Vieira.
\newblock A just-nonsolvable torsion-free group defined on the binary tree.
\newblock {\em J. Algebra}, 211:99--144, 1999.

\bibitem{gr_zu:lamp}
Rostislav~I. Grigorchuk and Andrzej {\.{Z}uk}.
\newblock The lamplighter group as a group generated by a $2$-state automaton
  and its spectrum.
\newblock {\em Geom. Dedicata}, 87(1--3):209--244, 2001.

\bibitem{bgr:spec}
Laurent Bartholdi and Rostislav~I. Grigorchuk.
\newblock On the spectrum of {H}ecke type operators related to some fractal
  groups.
\newblock {\em Proceedings of the Steklov Institute of Mathematics}, 231:5--45,
  2000.

\bibitem{cuntz_rep}
Ola Bratelli and Palle~E.~T. Jorgensen.
\newblock {\em Iterated function systems and Permutation representations of the
  {Cuntz} algebra}, volume 139 of {\em Memoirs of the American Mathematical
  Society}.
\newblock A.~M.~S., Providence, Rhode Island, 1999.

\bibitem{yacobson3}
M.~V. Yacobson.
\newblock On the question of topological classification of rational mappings of
  the {Riemann} sphere.
\newblock {\em Uspekhi Mat. Nauk}, 28(2):247--248, 1973.

\bibitem{zukgrigorchuk:3st}
Rostislav~I. Grigorchuk and Andrzej {\.{Z}uk}.
\newblock On a torsion-free weakly branch group defined by a three state
  automaton.
\newblock {\em Internat. J. Algebra Comput.}, 12(1):223--246, 2002.

\bibitem{DH:Thurston}
Adrien Douady and John~H. Hubbard.
\newblock A proof of {Thurston's} topological characterization of rational
  functions.
\newblock {\em Acta Math.}, 171(2):263--297, 1993.

\bibitem{neksid}
V.~Nekrashevych and S.~Sidki.
\newblock Automorphisms of the binary tree: state-closed subgroups and dynamics
  of $1/2$-endomorphisms.
\newblock In T.~W. {M\"uller} and H.~Helling, editors, {\em {Proceedings of the
  conference `Groups: combinatorial and geometric aspects', held in Bielefeld
  15-22 August 1999}}, LMS Lecture Notes Series, 2003.
\newblock to appear.

\bibitem{zukgrigorchuk:3stsp}
Rostislav~I. Grigorchuk and Andrzej {\.{Z}uk}.
\newblock On the spectrum of a torsion-free weakly branch group defined by a
  three state automaton.
\newblock In {\em Computational and Statistical Group Theory}, Contemp. Math.,
  Amer. Math. 2002.
\newblock to appear.

\end{thebibliography}

\end{document}